\documentclass{amsart}

\usepackage{multicol}
\usepackage{amsmath}
\usepackage[latin1]{inputenc}
\usepackage{amsfonts}
\usepackage{amssymb}
\usepackage{amsthm}
\usepackage{amscd}

\usepackage[latin1]{inputenc}

\usepackage{amsmath}
\usepackage{pb-diagram}
\usepackage[mathscr]{eucal}
\usepackage{amsfonts}

\def\sn{-\hspace{-.17cm}-}
\usepackage{amssymb}

\usepackage{amsthm}

\usepackage{graphics}

\newcommand{\mk}{\medskip}

\newcommand{\ZZ}{\mathbb{Z}}

\newcommand{\CC}{\mathbb{C}}

\newcommand{\NN}{\mathbb{N}}

\newcommand{\QQ}{\mathbb{Q}}

\newcommand{\Glie}{\mathfrak{g}}

\newcommand{\Yim}{\mathcal{Y}}

\newcommand{\Hlie}{\mathfrak{h}}

\newcommand{\demo}{\noindent {\it \small Proof:}\quad}

\newcommand{\U}{\mathcal{U}}

\newcommand{\Lo}{\mathcal{L}}

\newtheorem{thm}{Theorem}[section]

\newtheorem{defi}[thm]{Definition}

\newtheorem{cor}[thm]{Corollary}

\newtheorem{prop}[thm]{Proposition}

\newtheorem{lem}[thm]{Lemma}

\newtheorem{rem}[thm]{Remark}

\usepackage[all]{xy}
\title{Smallness problem for quantum affine algebras and quiver varieties}
\author{David Hernandez}
\address{CNRS - UMR 8100 : Laboratoire de Math\'ematiques de Versailles, 45 avenue des Etats-Unis , Bat. Fermat, 78035 VERSAILLES, 
FRANCE}
\email{hernandez@math.cnrs.fr}
\urladdr{http://www.math.uvsq.fr/\textasciitilde hernandez}

\begin{document}

\begin{abstract} The geometric small property (Borho-MacPherson \cite{bm}) of projective morphisms implies a description of their singularities in terms of intersection homology. In this paper we solve the smallness problem raised by Nakajima \cite{Nab, Nae} for certain resolutions of quiver varieties \cite{Nab} (analogs of the Springer resolution) : for Kirillov-Reshetikhin modules of simply-laced quantum affine algebras, we characterize explicitly the Drinfeld polynomials corresponding to the small resolutions. We use an elimination theorem for monomials of Frenkel-Reshetikhin $q$-characters that we establish for non necessarily simply-laced quantum affine algebras. We also refine results of \cite{her06} and extend the main result to general simply-laced quantum affinizations, in particular to quantum toro\"idal algebras (double affine quantum algebras).
\vskip 4.5mm

\noindent {\bf 2000 Mathematics Subject Classification:} Primary 17B37, Secondary 14L30, 81R50, 82B23, 17B67.
\\

\begin{center}PROBL\`EME DE PETITESSE POUR LES ALG\`EBRES AFFINES QUANTIQUES ET LES VARI\'ET\'ES CARQUOIS\end{center}

R\'ESUM\'E. La propri\'et\'e g\'eom\'etrique de petitesse (Borho-MacPherson \cite{bm}) des morphismes projectifs implique une description de leurs singularit\'es en termes d'homologie d'intersection. Dans cet article nous r\'esolvons le probl\`eme de petitesse pos\'e par Nakajima \cite{Nab, Nae} pour certaines résolutions de vari\'et\'es carquois \cite{Nab} (analogues de la r\'esolution de Springer) : pour les modules de Kirillov-Reshetikhin des alg\`ebres affines quantiques simplement lac\'ees, nous caract\'erisons explicitement les polyn\^omes de Drinfeld correspondant aux r\'esolutions petites. Nous utilisons un th\'eor\`eme d'\'elimination pour les mon\^omes des $q$-caract\`eres de Frenkel-Reshetikhin, que nous \'etablissons pour les alg\`ebres affines quantiques non n\'ecessairement simplement lac\'ees. Nous raffinons \'egalement des r\'esultats de \cite{her06} et \'etendons le r\'esultat principal aux affinis\'ees quantiques g\'en\'erales simplement lac\'ees, en particulier aux alg\`ebres toro\"idales quantiques (alg\`ebres quantiques doublement affines).

\end{abstract}

\maketitle

\tableofcontents

\section{Introduction} 

Borho and MacPherson introduced \cite[Section 1.1]{bm} remarkable geometric properties (smallness and semi-smallness) for a proper algebraic map $\pi : Z \rightarrow X$ where $Z, X$ are irreducible complex algebraic varieties : for a finite stratification of $X$ into irreducible smooth subvarieties, $\pi$ is said to be semi-small if the dimension of the inverse image of a point in a stratum is at most half the codimension of the stratum, and $\pi$ is said to be small if moreover the equality holds only if the stratum is dense. These properties does not depend of the stratification. 

\noindent This geometric situation is of particular interest as the Beilinson-Bernstein-Deligne-Gabber decomposition Theorem \cite{bbdg} is simplified \cite[Section 1.5]{bm} and provides an elegant description of the singularities of such maps in terms of intersection homology sheaves \cite{gm1, gm2}. A fundamental example of a semi-small morphism is given by the Springer resolution of the nilpotent cone of a complex simple Lie algebra, and the corresponding partial resolutions \cite{bm}. Nakajima \cite{nun, NaDuke} defined important and intensively studied varieties called quiver varieties which depend on a quiver $Q$. They come with a resolution which is semi-small \cite[Corollary 10.11]{NaDuke} for a finite Dynkin diagram (see \cite[Section 5.2]{NaCong}). 

\noindent The graded version of quiver varieties are also of particular importance, for example for their deep relations with representations of quantum affine algebras (see \cite{Nab}; the precise definition is reminded bellow). They also come with resolutions. A natural problem addressed in the present paper is to study the small property of these resolutions : in the present paper we address \cite[Conjecture 10.4]{Nab} (see also \cite{Nae}). Our study relies on the representation theory of quantum affine algebras. Let us also give the representation theoretical context for our study. 

In this paper $q\in\CC^*$ is fixed and is not a root of unity. Affine Kac-Moody algebras $\hat{\Glie}$ are infinite dimensional analogs of semi-simple Lie algebras $\Glie$, and have remarkable applications in several branches of mathematics and physics (see \cite{kac}). Their quantizations $\U_q(\hat{\Glie})$, called quantum affine algebras, have a very rich representation theory which has been intensively studied (see \cite{Cha2, dm} for references). In particular Drinfeld \cite{Dri2} discovered that they can also be realized as quantum affinization of usual quantum groups $\U_q(\Glie)$. By using this new realization, Chari-Pressley \cite{Cha2} classified their finite dimensional representations : they are parametrized by Drinfeld polynomials $(P_i(u))_{1\leq i\leq n}$ where $n$ is the rank of $\Glie$ and $P_i(u)\in\CC[u]$ satisfies $P_i(0) = 1$ .

\noindent A particular class of finite dimensional representations, called special modules, attracted much attention as Frenkel-Mukhin \cite{Fre2} proposed an algorithm which gives their $q$-character (analogs of usual characters adapted to the Drinfeld presentation of quantum affine algebras introduced by Frenkel-Reshetikhin \cite{Fre}). Let us give some examples : for $k > 0, i\in I, a\in\CC^*$, the Kirillov-Reshetikhin module $W_{k,a}^{(i)}$ is the simple module with Drinfeld polynomials 
\begin{equation*}
\begin{cases}
 P_j(u)=1 \text{ for }j\neq i, 
\\P_i(u)= (1-uaq_i^{k-1})(1-uaq_i^{k-3})\cdots (1-uaq_i^{1-k}).
\end{cases}
\end{equation*}
(The $q_i$ are certain power of $q$, see section \ref{un}). The $V_i(a)=W_{1,a}^{(i)}$ are called fundamental representations. The fundamental representations \cite{Fre2}, and the Kirillov-Reshetikhin modules \cite{Nad, her06} are special modules (this is the crucial point for the proof of the Kirillov-Reshetikhin conjecture). The corresponding standard module 
$$M(X_{k,a}^{(i)}) = V_i(aq_i^{1-k})\otimes V_i(aq_i^{3-k})\otimes\cdots \otimes V_i(aq_i^{k - 1})$$ 
is not special in general. 

\noindent The breakthrough geometric approach of Nakajima \cite{Naams, Nab} to $q$-characters of representations of simply-laced quantum affine algebras via (graded) quiver varieties led to remarkable advances in the description of finite dimensional representations : for example this approach provides an algorithm \cite{Nab} which computes the $q$-characters of {\it any} simple finite dimensional representations. Although in general the algorithm is very complicated, in some situations it provides a powerful tool to study these representations (for instance see \cite{Nad}). 

From the geometric point of view, the natural notion of small modules appeared in the following way : the small property of modules \cite{Nab} is the representation theoretical interpretation of the smallness of certain resolutions of (graded) quiver varieties mentioned above. 

\noindent A small module is special (but the converse is false in general). The representation theoretical interest of this notion is that all simple modules occurring in the Jordan-Hölder series of a small module are special, and so can be described by using the Frenkel-Mukhin algorithm.

\noindent A natural question is to characterize these small modules, and so the corresponding small resolutions. In particular, Nakajima (\cite[Conjecture 10.4]{Nab}, \cite{Nae}) raised the problem of characterizing the small standard modules corresponding to Kirillov-Reshetikhin modules.
\\

\noindent In this paper we solve this problem by giving explicitly the corresponding Drinfeld polynomials. 
\\

The crucial point for our proof is an elimination theorem for monomials of $q$-characters, that we establish by refining our results of \cite{her06}. Indeed it is easy to produce monomials that occur in a certain $q$-character (for example see remark \ref{process} bellow). But in general it is not clear if a given monomial {\it does not} occur in a $q$-character. The elimination theorem gives a criterion which implies that a monomial can be eliminated from the $q$-character of a simple module. Beyond the main result of the present paper (answer to the geometric smallness problem), we hope that this elimination theorem will be useful for other open problems in representation theory of quantum affine algebras. We already used it in a weak (non explicitly stated) form to prove the Kirillov-Reshetikhin conjecture \cite{her06}. Moreover it is used in \cite{miniaff} to study minimal affinizations of representations of quantum groups.

\noindent Let us state the main result of this paper. It can be stated in a simple compact way by using the following elementary definitions ($I=\{1,\cdots,n\}$ is the set of vertices of the Dynkin diagram of $\Glie$) : 

\begin{defi} A node $i\in \{1,\cdots ,n\}$ is said to be extremal (resp. special) if there is a unique $j\in I$ (resp. three distinct $j,k,l\in I$) such that $C_{i,j} < 0$ (resp. $C_{i,j}<0$, $C_{i,k}<0$ and $C_{i,l}<0$). 

For $i\in I$, we denote by $d_i$ the minimal $d\geq 1$ such that there are distinct $i_1,\cdots,i_d\in I$ satisfying $C_{i_j,i_{j+1}} < 0$ and $i_d$ is special (if there are no special vertices, we set $d_i = +\infty$ for all $i\in I$).\end{defi}
For example for $\Glie$ of type $A$, we have $d_i = +\infty$ for all $i\in I$. 

\noindent For illustration, examples are given on the following pictures :

Extremal node $i$ : {\large \vspace{-.15cm} 
$$\stackrel{i}{\circ}\hspace{-.18cm}\sn\hspace{-.18cm}\stackrel{j}{\circ} \hspace{-.18cm}\sn\hspace{-.18cm}{\circ} 
\dots $$}

Special node $i$ : {\large 
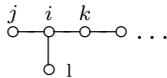
\vspace{-.15cm} 
$$\stackrel{j}{\circ}\hspace{-.18cm}\sn\hspace{-.18cm}\stackrel{i}{\circ} \hspace{-.18cm}\sn\hspace{-.18cm}\stackrel{k}{\circ} 
\hspace{-.18cm}\sn\hspace{-.18cm}\stackrel{}{\circ}\dots $$
\vspace{-.99cm}$$\hspace{-1.1cm}|$$
\vspace{-.99cm}$$\hspace{-.88cm}\circ\text{ \scriptsize l}$$}

Distance $d$ to a special node :{\large 
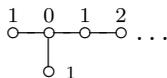
\vspace{-.15cm} 
$$\stackrel{1}{\circ}\hspace{-.18cm}\sn\hspace{-.18cm}\stackrel{0}{\circ} \hspace{-.18cm}\sn\hspace{-.18cm}\stackrel{1}{\circ} \hspace{-.18cm}\sn\hspace{-.18cm}\stackrel{2}{\circ}
\dots $$
\vspace{-1.06cm}$$\hspace{-1.1cm}|$$
\vspace{-.99cm}$$\hspace{-.83cm}\circ\text{ \scriptsize 1}$$}

\begin{thm}\label{mainres}[Smallness problem] Let $k > 0, i\in I, a\in\CC^*$. Then $M(X_{k,a}^{(i)})$ is small if and only if $k\leq 2$ or ($i$ is extremal and $k\leq d_i+1$).\end{thm} 

Remark : the condition is independent of the parameter $a\in\CC^*$.

In particular for $\Glie = sl_2$ or $\Glie = sl_3$, all $M(X_{k,a}^{(i)})$ are small (it proves the corresponding \cite[Conjecture 10.4]{Nab}). In general it gives an explicit criterion so that the smallness holds. On the geometric side, it characterizes the small resolutions mentioned above. 

\noindent Besides the statement of Theorem \ref{mainres} is also satisfied for all simply-laced quantum affinizations $\U_q(\hat{\Glie})$ ($\Glie$ is an arbitrary simply-laced Kac-Moody algebra, not necessarily semi-simple), in particular for quantum toro\"idal algebras (double affine quantum algebras).

The general idea of the proof is first to establish the result for type $A$ by using the elimination strategy of monomials explained above. We prove by induction on the highest weight that representations in a certain class are special. Then an argument allows to use the type $A$ to prove the result for general types.

Let us describe the organization of this paper. In section \ref{geom} we explain the geometric context of our results. In section \ref{un} we give some background on finite dimensional representations of quantum affine algebras and $q$-characters. In section \ref{secsmall} we recall from \cite{Nab} the definition of small modules and the geometric characterization (Theorem \ref{geomint}). We refine a Theorem of \cite{Nab} and give a more representation theoretical characterization (Theorem \ref{repthchar}). However this last result is not enough to prove Theorem \ref{mainres}, and technical work is needed in the next sections. The first point is the (representation theoretical) elimination Theorem (Theorem \ref{racourc}) which is proved in section \ref{seccan} : it gives a condition which implies that a monomial {\it does not} appear in the $q$-character of a simple module. Additional technical results are also proved in section \ref{seccan} : in particular the notion of thin modules (with $l$-weight spaces of dimension $1$) is introduced and studied. In section \ref{small}, we complete the proof of Theorem \ref{mainres} : first type $A$ is discussed, and then the general case is treated. The proof of the result for general simply-laced quantum affinizations is also discussed.

{\bf Acknowledgments :} The author is very grateful to Hiraku Nakajima for having attracted his attention to the smallness problem, and to Olivier Schiffmann for useful discussions.

\section{The geometric problem : small property and graded quiver varieties}\label{geom}

The geometric motivations and context of the results of the present paper have been explained at the beginning of the introduction. In this section we develop this discussion and define more precisely the involved geometric objects.

\subsection{Small property}

Let us recall the notion of semi-small and small morphism maps in the sense of Borho-MacPherson \cite{bm} for a proper algebraic map $\pi : Z \rightarrow X$ where $Z, X$ are irreducible complex algebraic varieties. 

\noindent We consider a finite stratification $X = \sqcup_i X_i$ into irreducible smooth subvarieties such that $\pi_{|\pi^{-1}(X_i)}$ is a topological fibration with base $X_i$ and fiber $\pi^{-1}(x_i)$ where $x_i\in X_i$ is a distinguished base point.

\begin{defi}\cite{bm} $\pi$ is said to be semi-small if for all $i$, 
$$2\text{dim}(\pi^{-1}(x_i))\leq \text{dim}(X) - \text{dim}(X_i).$$ 

$\pi$ is said to be small if $\pi$ is semi-small and if 
$$(2\text{dim}(\pi^{-1}(x_i)) = \text{dim}(X) - \text{dim}(X_i) \Rightarrow \text{dim}(X) = \text{dim}(X_i)).$$
In this case $X_i$ is said to be relevant.
\end{defi}
Note that stratification $X = \sqcup_i X_i$ exists (\cite{ha, t}) and that the property of being semi-small or small does not depend of the stratification.

\noindent When $\pi$ is projective and $Z$ is rationally smooth, this geometric situation is of particular interest as there is a very elegant description \cite[Section 1.5]{bm} of the singularities of such maps in terms of intersection homology sheaves \cite{gm1, gm2} : by using \cite[Section 1.7]{bm} the decomposition Theorem of Beilinson-Bernstein-Deligne-Gabber \cite{bbdg}, for $u\in X$, the cohomology groups $H^i(\pi^{-1}(u),\QQ)$ of the fiber $\pi^{-1}(u)$ are given by explicit formula involving the intersection homology of the closures $\overline{X_i}$ of strata such that $u\in\overline{X_i}$. The formula \cite[Section 1.5]{bm}can be expressed as a sum indexed by certain pairs $(X_i,\phi)$ where :

$X_i$ is a relevant stratum, 

$u\in\overline{X_i}$, 

$\phi$ is an irreducible representation of the fundamental group $\pi_1(X_i)$ of $X_i$, 

$\phi$ occurs in the decomposition of the representation of $\pi_1(X_i)$ on \\$H^{2\text{dim}(\pi^{-1}(x_i))}(\pi^{-1}(x_i),\QQ)$ by monodromy.

\noindent The case of small resolutions is remarkable, as the formula reduces to a single summand (and in this case the result is essentially given in \cite{gm2}). 

\noindent A fundamental example of semi-small morphism is given by the Springer resolution $T^*\mathcal{B}\rightarrow \mathcal{N}$ of the nilpotent cone of a complex simple Lie algebra, and the corresponding partial resolutions \cite{bm}. 

\noindent Nakajima \cite{nun, NaDuke} defined important and intensively studied varieties $\mathfrak{M}(v,w)$, $\mathfrak{M}_0(v,w)$ called quiver varieties which depend on a quiver $Q$ (see \cite{NaCong, os} for recent reviews). They come with a resolution 
$$\pi\colon \mathfrak{M}(v,w) \to \mathfrak{M}_0(v,w),$$ 
which gives an analog of the Springer resolution. It is proved in \cite[Corollary 10.11]{NaDuke} for $Q$ a finite Dynkin diagram type that $\pi$ is semi-small (see \cite[Section 5.2]{NaCong}). 

\noindent The graded version of quiver varieties $\mathfrak{M}^{\bullet}(V,W), \mathfrak{M}^{\bullet}_0(V,W)$ are also of particular importance, for example for their deep relations with representations of quantum affine algebras (see \cite{Nab}). 

Let us recall the definition of these varieties :

\subsection{(Graded) Quiver varieties} This section is essentially contained in \cite{Nab}. 

Fix a Dynkin diagram and an orientation on this diagram. Let $H$ be the set of oriented edges of the Dynkin diagram. For $h\in H$, $\operatorname{in}(h)$ (resp.\ $ \operatorname{out}(h)$)
is the incoming (resp.\ outgoing) vertex of $h$, and $\overline h$ is the same edge as $h$ with the
reverse orientation. We fix $q : H \rightarrow \{1, - 1\}$ such that $q(h) = - q(\overline{h})$ for any $h\in H$.

Let $V = \bigoplus_{i\in I,a\in\CC^*}V_{i,a}$ (resp. $W = \bigoplus_{i\in I,a\in\CC^*}W_{i,a}$) be a $I\times \CC^*$-graded vector spaces such that the $V_{i,a}$ (resp. $W_{i,a}$) are finite dimensional
and for at most finitely many $i\times a$, $V_{i,a}\neq 0$ (resp. $W_{i,a}\neq 0$). Consider for $n\in \ZZ$ :
$$  \operatorname{L}^\bullet(V, W)^{[n]} =
  \bigoplus_{i\in I, a\in\CC^*}
    \text{Hom} (V_{i,a}, W_{i,a q^{n}}),$$
$$  \operatorname{E}^\bullet(V, W)^{[n]} =
  \bigoplus_{h\in H, a\in\CC^*}
    \text{Hom} (V_{ \operatorname{out}(h),a}, W_{ \operatorname{in}(h),a q^{n}}),$$
$$  {\mathbf M}^\bullet(V, W) =
  \operatorname{E}^{\bullet}(V, V)^{[-1]} \oplus \operatorname{L}^{\bullet}(W, V)^{[-1]}
  \oplus \operatorname{L}^{\bullet}(V, W)^{[-1]}.$$
The above three components for an element of ${\mathbf M}^\bullet(V, W)$ are denoted by $B$,
$\alpha$, $\beta$ respectively,
the $\text{Hom} (V_{ \operatorname{out}(h),a},V_{ \operatorname{in}(h),a q^{-1}})$-component of 
$B$ is denoted by $B_{h,a}$ and we denote by $\alpha_{i,a}$,
$\beta_{i,a}$ the components of $\alpha$, $\beta$. Consider the map 
$$\mu\colon{\mathbf M}^\bullet(V,W)\to \operatorname{L}^{\bullet}(V,V)^{[-2]}$$ 
defined by
$$   \mu_{i,a}(B,\alpha,\beta)
   =  \sum_{ \operatorname{in}(h)=i}  q(h)
      B_{h,a q^{-1}} B_{\overline{h},a} +
   \alpha_{i,a q^{-1}}\beta_{i,a},$$
where $\mu_{i,a}$ is the $(i,a)$-component of $\mu$. We have an action of $G_V = \prod_{i,a} \operatorname{GL}(V_{i,a})$ on ${\mathbf M}^\bullet(V,W)$ defined by
$$  (B, \alpha, \beta) \mapsto g\cdot (B, \alpha, \beta)
  = (g_{ \operatorname{in}(h),a q^{-1}} B_{h,a} g_{ \operatorname{out}(h),a}^{-1},\,
  g_{i,a q^{-1}}\alpha_{i,a},\,
  \beta_{i,a} g_{i,a}^{-1} ).$$
The subvariety $\mu^{-1}(0)$ in ${\mathbf M}^{\bullet}(V,W)$ is stable under the action. 

Let us denote by $\mu^{-1}(0)^{\operatorname{s}}$ the set of stable points $(B, \alpha, \beta) \in \mu^{-1}(0)$, that is to say satisfying the condition : if an $I\times\CC^*$-graded subspace $S$ of $V$ is
$B$-invariant and contained in $\text{Ker}(\beta)$, then $S = 0$. The stability condition is invariant under the action of $G_V$, so we may say an orbit is stable or not. 

Consider the following quotient spaces of $\mu^{-1}(0)$:
\begin{equation*}
   \mathfrak{M}^{\bullet}_0(V,W) = \mu^{-1}(0)/\!\!/ G_V, \qquad
   \mathfrak{M}^{\bullet}(V,W) = \mu^{-1}(0)^{\operatorname{s}}/G_V.
\end{equation*}
Here $/\!\!/$ is the affine algebro-geometric quotient, the second one is the set-theoretical quotient. By \cite[3.18]{NaDuke}, there exists a natural projective morphism
\begin{equation*}
   \pi\colon \mathfrak{M}^{\bullet}(V,W) \to \mathfrak{M}^{\bullet}_0(V,W).
\end{equation*}
For $x\in \mu^{-1}(0)^{\operatorname{s}}$, $\pi(G_V.x)$ is the unique closed orbit contained in the closure of $G_V.x$.

\noindent Here $\mathfrak{M}^{\bullet}(V,W)$ is non singular and $\pi$ can be considered as an analog of the Springer resolution.

A natural problem addressed in the present paper is to study the small property of such resolutions $\pi$ : in the present paper we address \cite[Conjecture 10.4]{Nab} (see also \cite{Nae}). 

\noindent As our proof relies on the representation theory of quantum affine algebras, let us give some background about this subject :

\section{Quantum affine algebras and their representations}\label{un}

In this section we recall definitions and results about the representation theory of quantum affine algebras.

\subsection{Cartan matrix and quantized Cartan matrix} Let $C=(C_{i,j})_{1\leq i,j\leq n}$ be a Cartan matrix of 
finite type. We denote 
$I=\{1,\cdots,n\}$. $C$ is symmetrizable : there is a matrix $D=\text{diag}(r_1,\cdots,r_n)$ ($r_i\in\NN^*$)\label{ri} such 
that $B=DC$\label{symcar} is symmetric. In particular if $C$ is symmetric then $D=I_n$ (simply-laced case). 

\noindent We consider a realization $(\Hlie, \Pi, \Pi^{\vee})$ of $C$ (see \cite{bou, kac}): $\Hlie$ is a $n$ dimensional $\QQ$-vector space, $\Pi=\{\alpha_1,\cdots ,\alpha_n\}\subset \Hlie^*$ (set of the simple roots) and $\Pi^{\vee}=\{\alpha_1^{\vee},\cdots ,\alpha_n^{\vee}\}\subset \Hlie$ (set of simple coroots) are set such that for $1\leq i,j\leq n$, $\alpha_j(\alpha_i^{\vee})=C_{i,j}$.
Let $\Lambda_1,\cdots ,\Lambda_n\in\Hlie^*$ (resp. $\Lambda_1^{\vee},\cdots ,\Lambda_n^{\vee}\in\Hlie$) be the the fundamental weights (resp. coweights) : $\Lambda_i(\alpha_j^{\vee})=\alpha_i(\Lambda_j^{\vee})=\delta_{i,j}$ where $\delta_{i,j}$ is $1$ if $i=j$ and $0$ otherwise. Denote $P=\{\lambda \in\Hlie^*|\forall i\in I, \lambda(\alpha_i^{\vee})\in\ZZ\}$ the set of weights and $P^+=\{\lambda \in P|\forall i\in I, \lambda(\alpha_i^{\vee})\geq 0\}$ the set of dominant weights. For example we have $\alpha_1,\cdots ,\alpha_n\in P$ and $\Lambda_1,\cdots ,\Lambda_n\in P^+$. Denote $Q={\bigoplus}_{i\in I}\ZZ \alpha_i\subset P$ the root lattice and $Q^+={\sum}_{i\in I}\NN \alpha_i\subset Q$. For $\lambda,\mu\in \Hlie^*$, denote $\lambda \geq \mu$ if $\lambda-\mu\in Q^+$. Let $\nu:\Hlie^*\rightarrow \Hlie$ linear such that for all $i\in I$, we have $\nu(\alpha_i)=r_i\alpha_i^{\vee}$. For $\lambda,\mu\in\Hlie^*$, $\lambda(\nu(\mu))=\mu(\nu(\lambda))$. We use the enumeration of vertices of \cite{kac}.

\noindent We denote $q_i=q^{r_i}$ and for $l\in\ZZ, r\geq 0, m\geq m'\geq 0$ we define in $\ZZ[q^{\pm}]$ :
$$[l]_q=\frac{q^l-q^{-l}}{q-q^{-1}}\text{ , }[r]_q!=[r]_q[r-1]_q \cdots [1]_q\text{ ,
}\begin{bmatrix}m\\m'\end{bmatrix}_q=\frac{[m]_q!}{[m-m']_q![m']_q!}.$$

\noindent For $a,b\in\ZZ$, we denote $q^{a+b\ZZ}=\{q^{a+br}|r\in\ZZ\}$ and $q^{a+b\NN}=\{q^{a+br}|r\in\ZZ,r\geq 0\}$.

\noindent Let $C(z)$ be the quantized Cartan matrix defined by ($i\neq j\in I$): 
$$C_{i,i}(z)=z_i+z_i^{-1}\text{ ,
}C_{i,j}(z)=[C_{i,j}]_z.$$ 
$C(z)$ is invertible (see \cite{Fre}). We denote by $\tilde{C}(z)$ the inverse matrix of
$C(z)$ and by $D(z)$ the diagonal matrix such that for $i,j\in I$, $D_{i,j}(z)=\delta_{i,j}[r_i]_z$.

\subsection{Quantum algebras}\label{qkma}

\subsubsection{Quantum groups}

\begin{defi} The quantum group $\U_q(\Glie)$ is the $\CC$-algebra with generators $k_i^{\pm 1}$, $x_i^{\pm}$ ($i\in I$) and 
relations: 
$$k_ik_j=k_jk_i\text{ , } k_ix_j^{\pm}=q_i^{\pm C_{i,j}}x_j^{\pm}k_i,$$
$$[x_i^+,x_j^-]=\delta_{i,j}\frac{k_i-k_i^{-1}}{q_i-q_i^{-1}},$$
$${\sum}_{r=0\cdots  1-C_{i,j}}(-1)^r\begin{bmatrix}1-C_{i,j}\\r\end{bmatrix}_{q_i}(x_i^{\pm})^{1-C_{i,j}-r}x_j^{\pm}(x_i^{\pm})^r=0 \text{ (for $i\neq j$)}.$$
\end{defi}

\noindent This algebra was introduced independently by Drinfeld \cite{Dri1} and Jimbo \cite{jim}. It is remarkable that 
one can define a Hopf algebra structure on $\U_q(\Glie)$ by : 
$$\Delta(k_i)=k_i\otimes k_i,$$ 
$$\Delta(x_i^+)=x_i^+\otimes 1 + k_i\otimes x_i^+\text{ , }\Delta(x_i^-)=x_i^-\otimes 
k_i^{-1}+ 1\otimes x_i^-,$$ 
$$S(k_i)=k_i^{-1}\text{ , }S(x_i^+)=-x_i^+k_i^{-1}\text{ , }S(x_i^-)=-k_ix_i^-,$$ 
$$\epsilon(k_i)=1\text{ , }\epsilon(x_i^+)=\epsilon(x_i^-)=0.$$

\noindent Let $\U_q(\Hlie)$ be the commutative subalgebra of $\U_q(\Glie)$ generated by the $k_i^{\pm 1}$ ($i\in I$).

\noindent For $V$ a $\U_q(\Hlie)$-module and $\omega\in P$ we denote by $V_{\omega}$ the weight space of weight 
$\omega$ : 
$$V_{\omega}=\{v\in V|\forall i\in I, k_i.v=q_i^{\omega(\alpha_i^{\vee})}v\}.$$ 
In particular we have $x_i^{\pm}.V_{\omega}\subset V_{\omega \pm \alpha_i}$. We say that $V$ is $\U_q(\Hlie)$-diagonalizable if $V={\bigoplus}_{\omega\in P}V_{\omega}$ (in particular $V$ is of type $1$).

\subsubsection{Quantum loop algebras} We will use the second realization (Drinfeld realization) of the quantum loop 
algebra $\U_q(\Lo\Glie)$ (subquotient of the quantum affine algebra $\U_q(\hat{\Glie})$) :

\begin{defi}\label{defiaffi} $\U_q(\Lo\Glie)$ is the algebra with
generators $x_{i,r}^{\pm}$ ($i\in I, r\in\ZZ$), $k_i^{\pm 1}$ ($i\in I$), $h_{i,m}$ ($i\in I, m\in\ZZ-\{0\}$) and the
following relations ($i,j\in I, r,r'\in\ZZ, m,m' \in\ZZ-\{0\}$): 
$$\label{afcart}[k_i,k_j]=[k_i,h_{j,m}]=[h_{i,m},h_{j,m'}]=0,$$
$$k_ix_{j,r}^{\pm}=q_i^{\pm C_{i,j}}x_{j,r}^{\pm}k_i,$$
\begin{equation*}\label{aidecom}[h_{i,m},x_{j,r}^{\pm}]=\pm \frac{1}{m}[mB_{i,j}]_qx_{j,m+r}^{\pm},\end{equation*}
$$[x_{i,r}^+,x_{j,r'}^-]=\delta_{i,j}\frac{\phi^+_{i,r+r'}-\phi^-_{i,r+r'}}{q_i-q_i^{-1}},$$
$$x_{i,r+1}^{\pm}x_{j,r'}^{\pm}-q^{\pm B_{i,j}}x_{j,r'}^{\pm}x_{i,r+1}^{\pm}=q^{\pm B_{i,j}}x_{i,r}^{\pm}x_{j,r'+1}^{\pm}-x_{j,r'+1}^{\pm}x_{i,r}^{\pm},$$
$${\sum}_{\pi\in\Sigma_s}{\sum}_{k=0\cdots s}(-1)^k\begin{bmatrix}s\\k\end{bmatrix}_{q_i}x_{i,r_{\pi(1)}}^{\pm}\cdots x_{i,r_{\pi(k)}}^{\pm}x_{j,r'}^{\pm}x_{i,r_{\pi(k+1)}}^{\pm}\cdots x_{i,r_{\pi(s)}}^{\pm}=0,$$
where the last relation holds for all $i\neq j$, $s=1-C_{i,j}$, all sequences of integers $r_1,\cdots ,r_s$. $\Sigma_s$ is
the symmetric group on $s$ letters. For $i\in I$ and $m\in\ZZ$, $\phi_{i,m}^{\pm}\in \U_q(\Lo\Glie)$ is determined
by the formal power series in $\U_q(\Lo\Glie)[[z]]$ (resp. in $\U_q(\Lo\Glie)[[z^{-1}]]$): 
$${\sum}_{m\geq 0}\phi_{i,\pm m}^{\pm}z^{\pm m}=k_i^{\pm}\text{exp}(\pm(q-q^{-1}){\sum}_{m'\geq 1}h_{i,\pm m'}z^{\pm m'}),$$ 
and $\phi_{i,\mp m}^{\pm}=0$ for $m>0$. \end{defi}

\noindent $\U_q(\Lo\Glie)$ has a Hopf algebra structure (from the Hopf algebra structure of $\U_q(\hat{\Glie})$).

\noindent For $J\subset I$ we denote by $\U_q(\Lo\Glie_J)\subset\U_q(\Lo\Glie)$ the subalgebra generated by the $x_{i,m}^{\pm}$, $h_{i,m}$, $k_i^{\pm 1}$ for $i\in J$. $\U_q(\Lo\Glie_J)$ is a quantum loop algebra associated to the semi-simple Lie algebra $\Glie_J$ of Cartan matrix $(C_{i,j})_{i,j\in J}$. For example for $i\in I$, we denote $\U_q(\Lo\Glie_i)=\U_q(\Lo\Glie_{\{i\}})\simeq \U_{q_i}(\Lo sl_2)$.

\noindent The subalgebra of $\U_q(\Lo\Glie)$ generated by the $h_{i,m}, k_i^{\pm 1}$ (resp. by the $x_{i,r}^{\pm}$) is denoted by $\U_q(\Lo\Hlie)$ (resp. $\U_q(\Lo\Glie)^{\pm}$).

\subsection{Finite dimensional representations of quantum loop algebras} Denote by 
$\text{Rep}(\U_q(\Lo\Glie))$ the Grothendieck ring of (type $1$) finite 
dimensional representations of $\U_q(\Lo\Glie)$.

\subsubsection{Monomials and $q$-characters}\label{defimono} Let $V$ be a representation in $\text{Rep}(\U_q(\Lo\Glie))$. The subalgebra $\U_q(\Lo\Hlie)\subset\U_q(\Lo\Glie)$ is commutative, so we have :
$$V = {\bigoplus}_{\gamma=(\gamma_{i,\pm m}^{\pm})_{i\in I, m\geq 0}}V_{\gamma},$$
$$\text{where : }V_{\gamma}=\{v\in V|\exists p\geq 0, \forall i\in I, m\geq 0, (\phi_{i,\pm m}^{\pm}-\gamma_{i,\pm m}^{\pm})^p.v=0\}.$$
The $\gamma=(\gamma_{i,\pm m}^{\pm})_{i\in I, m\geq 0}$ are called $l$-weights (or pseudo-weights) and the $V_{\gamma}\neq \{0\}$ are called $l$-weight spaces (or pseudo-weight spaces) of $V$. One can prove \cite{Fre} that $\gamma$ is necessarily of the form :
\begin{equation*}{\sum}_{m\geq 0}\gamma_{i,\pm m}^{\pm}u^{\pm 
m}=q_i^{\text{deg}(Q_i)-\text{deg}(R_i)}\frac{Q_i(uq_i^{-1})R_i(uq_i)}{Q_i(uq_i)R_i(uq_i^{-1})},
\end{equation*}
where $Q_i,R_i\in \CC(u)$ satisfy $Q_i(0)=R_i(0)=1$. 

\noindent Consider the ring $\Yim = \ZZ[Y_{i,a}^{\pm}]_{i\in I, a\in\CC^*}$. The Frenkel-Reshetikhin $q$-characters morphism $\chi_q$ \cite{Fre} encodes the $l$-weights $\gamma$ (see also \cite{kn}). It is an injective ring morphism : 
$$\chi_q:\text{Rep}(\U_q(\Lo\Glie))\rightarrow \Yim$$
defined by
$$\chi_q(V) = {\sum}_{\gamma}\text{dim}(V_{\gamma})m_{\gamma},$$
where :
$$m_{\gamma}={\prod}_{i\in I, a\in\CC^*}Y_{i,a}^{q_{i,a}-r_{i,a}},$$
$$Q_i(u)={\prod}_{a\in\CC^*}(1-ua)^{q_{i,a}}\text{ , }R_i(u)={\prod}_{a\in\CC^*}(1-ua)^{r_{i,a}}.$$
The $m_{\gamma}$ are called monomials (they are analogs of weight). We denote by $A$ the set of monomials of $\ZZ[Y_{i,a}^{\pm}]_{i\in I, a\in\CC^*}$. For an $l$-weight $\gamma$, we denote $V_{\gamma}=V_{m_{\gamma}}$. We will also use the notation $i_r^p=Y_{i,q^r}^p$ for $i\in I$ and $r,p\in\ZZ$.

\noindent For $J\subset I$, $\chi_q^J$ is the morphism of $q$-characters for $\U_q(\Lo\Glie_J)\subset\U_q(\Lo\Glie)$. 
For a $m$ monomial we denote $u_{i,a}(m)\in\ZZ$ such that $m={\prod}_{i\in I, a\in\CC^*}Y_{i,a}^{u_{i,a}(m)}$. We also denote 
$$\omega(m)={\sum}_{i\in I, a\in\CC^*}u_{i,a}(m)\Lambda_i\text{ , }u_i(m)=\sum_{a\in\CC^*}u_{i,a}(m)\text{ , }u(m)=\sum_{i\in I}u_i(m).$$ 
$m$ is said to be 
$J$-dominant if for all $j\in J, a\in\CC^*$ we have $u_{j,a}(m)\geq 0$. An $I$-dominant monomials is said to be dominant.

\noindent Observe that $\chi_q, \chi_q^J$ can also be defined for finite dimensional $\U_q(\Lo \Hlie)$-modules in the same way.

\noindent In the following for $V$ a finite dimensional $\U_q(\Lo\Hlie)$-module, we denote by $\mathcal{M}(V)$ the set 
of monomials occurring in $\chi_q(V)$. 

\noindent For $i\in I, a\in\CC^*$, consider the analogs of simple roots for monomials : 
\begin{equation*}
\begin{split}
A_{i,a}=&Y_{i,aq_i^{-1}}Y_{i,aq_i}\prod_{\{j|C_{j,i}=-1\}}Y_{j,a}^{-1}
\\&\times\prod_{\{j|C_{j,i}=-2\}}Y_{j,aq^{-1}}^{-1}Y_{j,aq}^{-1}\prod_{\{j|C_{j,i}=-3\}}Y_{j,aq^2}^{-1}Y_{j,a}^{-1}Y_{j,aq^{-2}}^{-1}.
\end{split}
\end{equation*}

\noindent As the $A_{i,a}$ are algebraically independent \cite{Fre} (because $C(z)$ is invertible), for $M$ a product of $A_{i,a}^{\pm 1}$ we can define $v_{i,a}(M)\in\ZZ$ by $M={\prod}_{i\in I, a\in\CC^*}A_{i,a}^{-v_{i,a}(m)}$. We put $v_i(M)={\sum}_{a\in\CC^*}v_{i,a}(M)$ and $v(M)={\sum}_{i\in I}v_i(M)$.

\noindent For $\lambda\in Q^+$ we set $v(\lambda)=-\lambda(\Lambda_1^{\vee}+\cdots +\Lambda_n^{\vee})$. For $M$ a 
product of $A_{i,a}^{\pm 1}$, we have $v(M)=v(\omega(\lambda))$.

\noindent For $m,m'$ two monomials, we write $m'\leq m$ if $m'm^{-1}$ is a product of $A_{i,a}^{-1}$. 

\begin{defi}\label{monomrn}\cite{Fre2} A monomial $m\in A-\{1\}$ is said to be right-negative if for all $a\in\CC^*$, 
for $L=\text{max}\{l\in\ZZ| \exists i\in I, u_{i,aq^l}(m)\neq 0\}$, we have $\forall j\in I$, $u_{j,aq^L}(m)\leq 0$.\end{defi}

\noindent Observe that a right-negative monomial is not dominant. 

\begin{lem}\cite{Fre2}\label{rn} 1) For $i\in I, a\in\CC^*$, $A_{i,a}^{-1}$ is right-negative.

2) A product of right-negative monomials is right-negative.

3) If $m$ is right-negative, then $m'\leq m$ implies that $m'$ is right-negative.\end{lem}

\noindent For $J\subset I$ and $Z\in \Yim$, we denote $Z^{\rightarrow J}$ the element of $\Yim$ obtained from $Z$ by putting $Y_{j,a}^{\pm 1}=1$ for $j\notin J$.

\subsubsection{$l$-highest weight representations}

\noindent The irreducible finite dimensional $\U_q(\Lo\Glie)$-modules have been classified by Chari-Pressley. They are parameterized by dominant monomials : 

\begin{defi}\label{lhigh} A $\U_q(\Lo\Glie)$-module $V$ is said to be of $l$-highest weight $m\in A$ if there is $v\in V_m$ such that $V=\U_q(\Lo\Glie)^-.v$ and $\forall i\in I, r\in\ZZ, x_{i,r}^+.v=0$.\end{defi}

\noindent For $m\in A$, there is a unique simple module $L(m)$ of $l$-highest weight $m$. 

\begin{thm}\cite[Theorem 12.2.6]{Cha2} The dimension of $L(m)$ is finite if and only if $m$ is dominant.\end{thm}

For $i\in I$, $a\in\CC^*$, $k\geq 0$ we denote $X_{k,a}^{(i)}={\prod}_{k'\in\{1,\cdots ,k\}}Y_{i,aq_i^{k-2k'+1}}$. 

\begin{defi} The Kirillov-Reshetikhin modules are the $W_{k,a}^{(i)} = L(X_{k,a}^{(i)})$.\end{defi}

For $i\in I$ and $a\in\CC^*$, $W_{1,a}^{(i)}$ is called a fundamental representation and is denoted by $V_i(a)$ (in the case $\Glie = sl_2$ we simply write $W_{k,a}$ and $V(a)$). 

For $m\in\ZZ[Y_{i,a}]_{i\in I, a\in\CC^*}$ a dominant monomial, the standard module $M(m)$ is defined \cite{Naams, vv} as the tensor product : 
$$M(m)=\bigotimes_{a\in(\CC^*/q^{\ZZ})}(\cdots\otimes(\bigotimes_{i\in I}V_i(aq)^{\otimes u_{i,aq}(m)})\otimes(\bigotimes_{i\in I}V_i(aq^2)^{\otimes u_{i,aq^2}(m)})\otimes\cdots).$$
It is well-defined as for $i,j\in I$ and $a\in\CC^*$ we have $V_i(a)\otimes V_j(a)\simeq V_j(a)\otimes V_i(a)$ and for $a'\notin aq^{\ZZ}$, we have $V_i(a)\otimes V_j(a')\simeq V_j(a')\otimes V_j(a)$. Observe that fundamental representations are particular cases of standard modules.

Let $\Glie = sl_2$. The monomials $m_1=X_{k_1,a_1}$, $m_2=X_{k_2,a_2}$ are said to be in special position if the monomial $m_3={\prod}_{a\in\CC^*}Y_a^{\text{max}(u_{a}(m_1),u_{a}(m_2))}$ is of the form $m_3=X_{k_3,a_3}$ and $m_3\neq m_1, m_3\neq m_2$. A normal writing of a dominant monomial $m$ is a product decomposition $m={\prod}_{i=1,\cdots,L}X_{k_l,a_l}$ such that for $l\neq l'$, $X_{k_l,a_l}$, $X_{k_{l'},a_{l'}}$ are not in special position. Any dominant monomial has a unique normal writing up to permuting the monomials (see \cite[Section 12.2]{Cha2}). It follows from the study of the representations of $\U_q(\Lo sl_2)$ in \cite{Cha0, Cha, Fre} that :

\begin{prop}\label{aidesldeux} Suppose that $\Glie=sl_2$.

(1) $W_{k,a}$ is of dimension $k+1$ and :
$$\chi_q(W_{k,a})=X_{k,a}(1+A_{aq^{k}}^{-1}(1+A_{aq^{k-2}}^{-1}(1+\cdots (1+A_{aq^{2-k}}^{-1}))\cdots).$$

(2) $V(aq^{1-k})\otimes V(aq^{3-k})\otimes \cdots\otimes V(aq^{k-1})$ is of $q$-character :
$$X_{k,a}(1+A_{aq^{k}}^{-1})(1+A_{aq^{k-2}}^{-1})\cdots (1+A_{aq^{2-k}}^{-1}).$$
In particular all $l$-weight spaces of the tensor product are of dimension $1$.

(3) For $m$ a dominant monomial and $m=X_{k_1,a_1}\cdots X_{k_l,a_l}$ a normal writing we have :
$$L(m)\simeq W_{k_1,a_1}\otimes \cdots\otimes W_{k_l,a_l}.$$
\end{prop}

\subsubsection{Special modules and complementary reminders}\label{compred}

Let us consider analogs of cones of weights (for example used to define category $\mathcal{O}$ for affine Kac-Moody algebras) adapted to monomials :

\begin{defi} For $m\in A$, $D(m)$ is the set of monomials $m'\in A$ such
that there are $m_0=m, m_1, \cdots, m_N=m'\in A$ satisfying for all $j\in\{1,\cdots ,N\}$ : 
\begin{enumerate}
\item $m_j=m_{j-1}A_{i_j,a_1q_{i_j}}^{-1}\cdots A_{i_j,a_{r_j}q_{i_j}}^{-1}$ where $i_j\in I$, $r_j\geq 1$ and $a_1,\cdots, a_{r_j}\in\CC^*$,

\item for $1\leq r\leq r_j$, $u_{i_j,a_r}(m_{j-1})\geq
|\{r'\in\{1,\cdots ,r_j\}|a_{r'}=a_r\}|$ where $r_j,i_j,a_r$ are as in condition (1).
\end{enumerate}\end{defi}

\noindent The motivation for this definition comes from the two simple facts :

for all $m'\in D(m)$, $m'\leq m$,

if $m'\in D(m)$, then $(D(m')\subset D(m))$,

\noindent and from the following result which gives a strong condition for a monomial to appear in a $q$-character :

\begin{thm}\label{inducb}\cite[Theorem 5.21]{her07} For $V$ a finite dimensional $l$-highest weight module of highest monomial $m$, we have $\mathcal{M}(V)\subset D(m)$.\end{thm}

\noindent In particular for all $m'\in\mathcal{M}(V)$, we have $m'\leq m$ and the $v_{i,a}(m'm^{-1}), v(m'm^{-1})\geq 0$ are well-defined. As a direct consequence of Theorem \ref{inducb}, we also have : 

\begin{lem}\label{fundless} For $i\in I, a\in\CC^*$, we have $(\chi_q(V_i(a))-Y_{i,a})\in\ZZ[Y_{j,aq^l}^{\pm}]_{j\in I, l > 0}$.\end{lem}

\noindent This last result was first proved in \cite[Lemma 6.1, Remark 6.2]{Fre2}.

\noindent The notion of special module was introduced in \cite{Nab} :

\begin{defi} A $\U_q(\Lo\Glie)$-module is said to be special if his $q$-character has a unique dominant monomial.\end{defi}

\noindent This notion is of particular importance because an algorithm of Frenkel-Mukhin \cite{Fre2} gives the $q$-character of special modules. Observe that a special module is a simple $l$-highest weight module (as each simple module occurring in the Jordan-Hölder series of a representation contributes with at least one dominant monomial in the $q$-character). But in general all simple $l$-highest weight module are not special. 

\noindent The following result was proved in \cite{Nab, Nad} for simply laced types, and in full generality in \cite{her06} (see \cite{Fre2} for previous results). It gives a remarkable example of a family of special modules and is the crucial point for the proof of the Kirillov-Reshetikhin conjecture :

\begin{thm}\cite[Theorem 4.1, Lemma 4.4]{her06}\label{formerkr} The Kirillov-Reshetikhin modules are special. Moreover for $m\in \mathcal{M}(W_{k,a}^{(i)})$, $m\neq X_{k,a}^{(i)}$ implies $m\leq X_{k,a}^{(i)}A_{i,aq_i^k}^{-1}$.\end{thm}

Now let us recall a decomposition result of $q$-characters relatively to sub-Dynkin diagrams corresponding to $J\subset I$ (Proposition \ref{jdecomp}). This is the analog at the level of $q$-character of the decomposition of a simple representation in simple representations for the subalgebra $\U_q(\Lo\Glie_J)$. This result will be intensively used in the following. 

\noindent Define $$\mu_J^I:\ZZ[(A_{j,a}^{\pm})^{\rightarrow (J)}]_{j\in J,a\in\CC^*}\rightarrow \ZZ[A_{j,a}^{\pm}]_{j\in J,a\in\CC^*},$$ 
the ring morphism such that $\mu_J^I((A_{j,a}^{\pm})^{\rightarrow(J)})=A_{j,a}^{\pm}$. For $m$ $J$-dominant, denote by $L^J(m^{\rightarrow(J)})$ the simple $\U_q(\Lo\Glie_J)$-module of $l$-highest weight $m^{\rightarrow(J)}$. Define :
$$L_J(m)=m \mu_J^I((m^{\rightarrow(J)})^{-1}\chi_q^J(L^J(m^{\rightarrow(J)}))).$$
(Observe that from Proposition \ref{aidesldeux}, we have explicit formulas for the $L_{\{i\}}(m)$ for $i\in I$.)

\begin{prop}\label{jdecomp}\cite[Proposition 3.1]{her05} For a representation $V\in\text{Rep}(\U_q(\Lo \Glie))$ and $J\subset I$, there is unique
decomposition in a finite sum : 
\begin{equation}\label{formdecomp}\chi_q(V)=\sum_{m'\text{ $J$-dominant}}\lambda_J(m')L_J(m').\end{equation}
Moreover for all $m'$ $J$-dominant we have $\lambda_J(m')\geq 0$.\end{prop}
(In \cite{her05} the $\lambda_J(m')\geq 0$ were assumed, but the proof of the uniqueness does not depend on it.)

\noindent As a consequence :

\begin{cor}\label{hproc} Let $m$ be a dominant monomial and $m'$ such that

(i) $m'\in\mathcal{M}(L(m))$,

(ii) $m'$ is $J$-dominant monomial, 

(iii) there are no $m'' > m'$ satisfying $m''\in\mathcal{M}(m)$ and $m'$ appears in $L_J(m'')$. 

\noindent Then the monomials of $L_J(m')$ are in $\mathcal{M}(L(m))$.
\end{cor}

\demo From the last condition $L_J(m')$ occurs in the decomposition of Proposition \ref{jdecomp}. As the coefficients in this decomposition are positive, all monomials on $L_J(m')$ occur in $\chi_q(L(m))$.\qed

\begin{rem}\label{process} In Corollary \ref{hproc}, we can start with $m' = m$, and then we use for $m'$ monomials in $L_J(m)$, and so on. This process gives inductively from $m$ a set of monomial occurring in $\chi_q(L(m))$.\end{rem}

\section{Representation theoretical interpretation of the small property}\label{secsmall}

In this section $\Glie$ is simply laced. 

Originally the notion of small modules was given in terms of $q,t$-characters \cite{Nab}. We recall this definition and the relation \cite{Nab} with the geometric small property of Section \ref{geom} (Theorem \ref{geomint}). 

Although the representation theoretical meaning of $q,t$-character is not totally understood (see \cite[Conjecture 3.1.1]{Naex}), the notion of small modules can be purely algebraically formulated : we give an additional representation theoretical interpretation of the notion (Theorem \ref{repthchar}) by refining a proof of \cite{Nab} (this provides an additional algebraic motivation for the study of the small modules).

We also comment the main result of the present paper (Theorem \ref{mainres}).

\subsection{Definition of small modules and $q,t$-characters} 

The notion of small modules is related to the notion of $q,t$-characters defined in \cite{Naex, Nab}. There are $t$-deformations of $q$-characters which can be purely algebraically defined (see \cite{her02} for non-simply laced cases with a different approach including a purely algebraic proof of the existence). They are a very powerful tool as Nakajima proved they provide an algorithm which allows to compute the $q$-character of any simple representation.

Consider the commutative ring $\hat{\Yim}_t=\ZZ[V_{i,a},W_{i,a},t^{\pm}]_{i\in I, a\in\CC^*}$. A monomial of $\hat{\Yim}_t$ is a product of $V_{i,a}$, $W_{i,a}$. One says $m'\leq m$ if $m'm^{-1}$ is a product of $V_{i,a}$. The $q,t$-characters map $\chi_{q,t} : \text{Rep}(\U_q(\Lo\Glie))\rightarrow \hat{\Yim}_t$ is a $\ZZ$-linear map defined by three axioms in \cite{Nab} : 

1) the data of the image of $\chi_{q,t}$, 

2) a compatibility property of the tensor product with a certain twisted product on $\hat{\Yim}_t$, 

3) for $m\in\ZZ[Y_{i,a}]_{i\in I, a\in\CC^*}$ a dominant monomial of $\Yim$, the relation :
$$\chi_{q,t}(M(m))\in M_0 + \sum_{m' < M_0}\ZZ[t^{\pm}] m'\text{ where }M_0=\prod_{i\in I,a\in\CC^*}W_{i,a}^{u_{i,a}(m)}.$$
(Only the last axiom will be explicitly used in the following, and so we refer to \cite{Nab} for the details of the first two axioms).

Let $m$ be a monomial of $\hat{\Yim}_t$. For $i\in I$, $a\in\CC^*$, one defines $w_{i,a}(m),v_{i,a}(m)\geq 0$ by $m=\prod_{i\in I, a\in\CC^*}W_{i,a}^{w_{i,a}(m)}V_{i,a}^{v_{i,a}(m)}$, and :
$$u_{i,a}(m)=w_{i,a}(m)-v_{i,aq^{-1}}(m)-v_{i,aq}(m)+\sum_{j\in I}C_{i,j}v_{j,a}(m),$$
$$d(m)=\sum_{i\in I,a\in\CC^*}(v_{i,aq}(m)u_{i,a}(m)+w_{i,aq}(m)v_{i,a}(m)).$$
We define a $\ZZ$-linear map $\hat{\Pi}:\hat{\Yim}_t\rightarrow\Yim$ by ($m$ is a monomial) :
$$\hat{\Pi}(m)=\prod_{i\in I,a\in\CC^*}Y_{i,a}^{u_{i,a}(m)}\text{ , }\hat{\Pi}(t)=1.$$ 
It is clear that $\hat{\Pi}$ is a ring morphism. 

\noindent A monomial $m$ of $\hat{\Yim}_t$ is said to be dominant if $\hat{\Pi}(m)$ is dominant. For $m$ a dominant monomial of $\hat{\Yim}_t$, one defines $M_t(m)\in\hat{\Yim}_t$ by :
$$M_t(m)=t^{d(m)}m(\prod_{i\in I, a\in\CC^*}W_{i,a}^{-u_{i,a}(m)})\chi_{q,t}(M(\prod_{i\in I, a\in\CC^*}Y_{i,a}^{u_{i,a}(m)}))\in\hat{\Yim}_t.$$

\begin{lem} For $m$ a dominant monomial, we have $\hat{\Pi}(M_t(m))=\chi_q(M(\hat{\Pi}(m)))$.\end{lem}

\demo From the defining axioms of $q,t$-characters, the evaluation at $t=1$ give $q$-characters \cite{Nab}, that is to say :
$$\hat{\Pi}(\chi_{q,t}(M(\prod_{i\in I, a\in\CC^*}Y_{i,a}^{u_{i,a}(m)}))) = \chi_q(M(\prod_{i\in I, a\in\CC^*}Y_{i,a}^{u_{i,a}(m)})) = \chi_q(M(\hat{\Pi}(m))).$$
As $\hat{\Pi}(t^{d(m)}m(\prod_{i\in I, a\in\CC^*}W_{i,a}^{-u_{i,a}(m)})) = 1$, the result is clear.
\qed

For $m$ a dominant monomial and $m'\leq m$ a monomial, $c_{m,m'}(t)\in \ZZ[t^{\pm}]$ is defined by :
$$M_t(m)=\sum_{m'\leq m}c_{m,m'}(t)t^{d(m')}m'.$$

\begin{defi}\label{original}\cite{Nab} Let $m$ be a dominant monomial of $\Yim$. The standard module $M(m)$ is said to be small if for all dominant monomials $m',m''\leq m$, we have $c_{m',m''}(t)\in t^{-1}\ZZ[t^{-1}]$.\end{defi}

Remark : Observe that in general there is no hope to have $c_{m',m''}(t)\in t^{-1}\ZZ[t^{-1}]$ without assuming that $m''$ is dominant. For example for $\Glie = sl_2$, we have 
$$\chi_{q,t}(M(W_a))=W_a+W_aV_{aq}\text{ , }d(W_aV_{aq})=0\text{ , }c_{W_a,W_aV_{aq}}=1\notin t^{-1}\ZZ[t^{-1}].$$ However $M(W_a)$ is small (see Proposition \ref{fundsmall} bellow).

\subsection{Geometric characterization}

The motivation for this Definition \ref{original} comes from geometry \cite{Nab} and from the relation to the small property of Section \ref{geom} :

Consider the monomials $m_W$, $m_V\in\hat{\Yim}_t$ defined by
$$M_W = \prod_{i\in I,a\in\CC^*} W_{i,a}^{\text{dim}(W_{i,a})}\text{ , }m_V = \prod_{i\in I,a\in\CC^*} V_{i,a}^{\text{dim}(V_{i,a})}.
$$
As a consequence of the geometric construction of representations of quantum affine algebras, we have the following geometric characterization of small standard modules (see \cite[Remark 10.2]{Nab}) :

\begin{thm}\label{geomint}\cite{Nab} Let $m$ be a dominant monomial of $\hat{\Yim}_t$ and $W = \bigoplus_{i\in I, a\in\CC^*}W_{i,a}$ be the graded space satisfying $\text{dim}(W_{i,a}) = u_{i,a}(m)$. The standard module $M(m)$ is small if and only if for all $V$ such that $M_W m_V$ is dominant, the resolution $\pi\colon \mathfrak{M}^{\bullet}(V,W) \to \mathfrak{M}^{\bullet}_0(V,W)$ is small.\end{thm}

\subsection{Representation theoretical characterization} Let us give another characterization of small modules.

Consider the $\ZZ$-linear involution of $\hat{\Yim}_t$ defined by $\overline{m}=t^{2d(m)}m$, $\overline{t}=t^{-1}$. Observe that for $m$ a monomial of $\hat{\Yim}_t$, $t^{d(m)}m$ is invariant by the involution.

In \cite{Nab} Nakajima constructed a family $\mathcal{L}(m)\in\hat{\Yim}_t$, indexed by the set of dominant monomial $m$ of $\hat{\Yim}_t$, characterized by the properties :

i) $\overline{\mathcal{L}(m)}=\mathcal{L}(m)$,

ii) $\mathcal{L}(m)\in M_t(m)+\sum_{\{m'\text{ dominant}|m' < m\}} t^{-1}\ZZ[t^{-1}]M_t(m')$.

\noindent They are analogs of canonical bases in $\hat{\Yim}$ for the bar involution, and the transition coefficient to the basis $(M_t(m))_m$ are analogs of Kazhdan-Lusztig polynomials.

Nakajima proved \cite{Nab} the following deep result : 

\begin{thm}\label{char}\cite{Nab} For all $m$ dominant monomial of $\hat{\Yim}_t$, we have $$\hat{\Pi}(\mathcal{L}(m))=\chi_q(L(\hat{\Pi}(m))).$$
\end{thm}

In particular this provides an algorithm to compute the $q$-characters of simple modules. It is very complicated in general, and it is difficult to get explicit formulas from it, but it provides applications in situations where the algorithm can be simplified (for example see \cite{Nad}).

As a consequence of this result, we have :

\begin{thm}\label{qtsimp}\cite{Nab} Let $m$ be a dominant monomial monomial of $\Yim$. If $M(m)$ is small, then for all dominant monomial $m'\leq m$, $L(m')$ is special.\end{thm}

In fact the converse is true by using the following two results :

\begin{thm}\label{posstandard}\cite[Theorem 3.5 (6)]{Nab} For all dominant monomial $m$ of $\hat{\Yim}_t$, the coefficient of a monomial occurring in $M_t(m)$ is a Laurent polynomial with nonnegative coefficients.\end{thm}

\begin{lem}\label{fini} For $M$ a dominant monomial of $\hat{\Yim}_t$, the set 
$$\{M'| M'\leq M \text{ and $M'$ is dominant}\}$$ 
is finite. \end{lem}

\demo We can suppose that $M\in\ZZ[Y_{i,aq^r}]_{i\in I,r\in\ZZ}$ where $a\in\CC^*$. Let $K = \text{max}\{r\in\ZZ|\exists i\in I,u_{i,aq^r}(M)\neq 0\}$. For $M' < M$, $M'M^{-1}$ is right-negative so $M'$ dominant implies $\sum_{i\in I,r\geq K}v_{i,aq^r} (M'M^{-1})= 0$. It is proved in \cite[Lemma 3.14]{her02} that the set 
$$\{M' = MA_{i_1,aq^{l_1}}^{-1}\cdots A_{i_R,aq^{l_R}}^{-1}|R\geq 0\text{ , }l_1,\cdots,l_R\leq K, M' \text{ is dominant}\}$$ 
is finite, and so we can conclude (note that in \cite{her02}, $l_1,\cdots,l_R\leq K$ is replaced by $l_1,\cdots,l_R\geq K$, but the proof is the same).\qed

By using a slight modification of the proof of Theorem \ref{qtsimp} in \cite{Nab}, we get the following characterization :

\begin{thm}\label{repthchar} Let $m$ be a dominant monomial of $\Yim$. $M(m)$ is small if and only if for all dominant monomial $m'\leq m$, $L(m')$ is special.\end{thm}

Observe that it is a purely representation theoretical characterization of small modules involving $q$-characters, without $q,t$-characters. This provides an additional algebraic motivation for the study of the small modules : all simple module which could appear in the "cone of monomial" of a small module are special, and so can be described by using the Frenkel-Mukhin algorithm.

\demo The only if part is the statement of Theorem \ref{char}. Let us prove the if part.

\noindent For $m,m'$ dominant monomials of $\hat{\Yim_t}$, we consider $Z_{m,m'}(t)\in\ZZ[t^{\pm}]$ defined by 
$$M_t(m)=\sum_{m'\text{ dominant}}Z_{m,m'}(t)\mathcal{L}(m').$$
By definition of $\mathcal{L}(m')$ we have $Z_{m,m}(t)=1$ and $Z_{m,m'}(t) \in t^{-1}\ZZ[t^{-1}]$ for $m' < m$. If $m' \nleq m$, we have $Z_{m,m'}(t) = 0$.

As $M = m (\prod_{i\in I, a\in\CC^*}W_{i,a}^{u_{i,a}(m)})^{-1}$ satisfies $u_{i,a}(M) = 0$ for any $i\in I, a\in\CC^*$, we can suppose that $m=\prod_{i\in I, a\in\CC^*}W_{i,a}^{u_{i,a}(m)}$. Suppose that for all dominant monomial $m'\leq m$, $L(\hat{\Pi}(m'))$ is special. From Lemma \ref{fini}, there is a finite number of dominant monomial $m'\leq m$. Choose a numbering $m_1,m_2,\cdots, m_N = m$ of these monomials such that $m_r < m_{r'}$ implies $r < r'$. Denote by $Z_{m_r,m_{r'}}'(t)$ the coefficients of $m_r'$ in $\mathcal{L}(m_r)$. As $L(\hat{\Pi}(m_r))$ is special, it follows from Theorem \ref{qtsimp} that $Z_{m_r,m_{r'}}'(1)=0$. Let us prove by induction $r$ that 
$$\forall r'' < r'\leq  r\text{ , }c_{m_{r'},m_{r''}}(t)\in t^{-1}\ZZ[t^{-1}]\text{ and }Z_{m_{r'},m_{r''}}'(t)=0.$$ 
For $r=1$, we have $\mathcal{L}(m_1)=M_t(m_1)$. Now we consider $r > 1$. By the induction hypothesis, for all $r' < r$, $\mathcal{L}(m_{r'})$ has no dominant monomial except $m_{r'}$. So for $r'' < r$, 
$$c_{m_r,m_{r''}}(t) = Z_{m_r,m_{r''}}(t) + Z_{m_r,m_{r''}}'(t).$$ 
From Theorem \ref{posstandard}, we have 
$$c_{m_r,m_{r''}}(t) = tP^+(t)+\alpha + t^{-1}P^-(t^{-1})$$ 
where $\alpha\geq 0$, $P^+,P^-\in\NN[t]$. As $Z_{m_r,m_{r''}}(t)\in t^{-1}\ZZ[t^{-1}]$ and $Z_{m_r,m_{r''}}'(t)=Z_{m_r,m_{r''}}'(t^{-1})$, we have 
$$Z_{m_r,m_{r''}}' (t) = tP^+(t) +\alpha +t^{-1}P^+(t^{-1}).$$ 
So $Z_{m_r,m_{r''}}'(t)$ has positive coefficients, so $Z_{m_r,m_{r''}}'(1)=0$ implies $Z_{m_r,m_{r''}}'(t)=0$. So $c_{m_r,m_{r''}}(t) = Z_{m_r,m_{r''}}(t)\in t^{-1}\ZZ[t^{-1}]$. As a conclusion, $M(m)$ is small.\qed

\subsection{Main result}

A natural question is to characterize small modules and so the corresponding small resolutions. In particular, Nakajima \cite[Conjecture 10.4]{Nab}, \cite{Nae} raised the problem of characterizing the Drinfeld polynomials of small standard modules corresponding to Kirillov-Reshetikhin modules.

The main result of this paper is an explicit answer to this question (Theorem \ref{mainres}). First let us note in general the standard modules corresponding to Kirillov-Reshetikhin modules are not necessarily small :

\begin{rem}\label{counter} Let $\Glie = sl_4$ and $m=Y_{2,1}Y_{2,q^2}Y_{2,q^4}$. Consider $m' = mA_{2,q}^{-1} = Y_{1,q}Y_{3,q}Y_{2,q^4}$. Then by using the process described in remark \ref{process}, the monomials $Y_{1,q^3}^{-1}Y_{3,q^3}^{-1}Y_{2,q^2}^2Y_{2,q^4}=m'A_{1,q^2}^{-1}A_{3,q^2}^{-1}$ and $Y_{2,q^2}=m'A_{1,q^2}^{-1}A_{3,q^2}^{-1}A_{2,q^3}^{-1}$ occur in $\chi_q(L(m'))$ and $L(m')$ is not special. 
So $M(m)$ is not small.\end{rem}

A crucial step for the proof of Theorem \ref{mainres} is the elimination theorem proved in the next section.

\section{Elimination theorem and preliminary results}\label{seccan}

In this section $\Glie$ is an arbitrary semi-simple Lie algebra. We prove several preliminary results so that we can prove Theorem \ref{mainres} in the last section of the paper.

\subsection{Elimination Theorem}

We have seen a (combinatorial) procedure which allows to produce monomials occurring in a $q$-character (remark \ref{process}). We first prove in this section a (representation theoretical) theorem (Theorem \ref{racourc}) which gives a criterion so that a monomial $m'$ does not occur in the $q$-character of a simple modules $L(m)$. 

\noindent This theorem is used in \cite{miniaff} to study minimal affinizations of representations of quantum groups.

\subsubsection{Statement}

\begin{thm}\label{racourc} Let $V = L(m)$ be a $\U_q(\Lo \Glie)$-module simple module. Let $m' < m$ and $i\in I$ satisfying the following conditions :

(i) there is a unique $i$-dominant $M\in(\mathcal{M}(V)\cap m'\ZZ[A_{i,a}]_{a\in\CC^*})-\{m'\}$ and its coefficient is $1$,

(ii) $\sum_{r\in\ZZ} x_{i,r}^+ (V_M)=\{0\}$,

(iii) $m'$ is not a monomial of $L_i(M)$,

(iv) if $m''\in\mathcal{M}(\U_q(\Lo\Glie_i).V_M)$ is $i$-dominant, then $v(m''m^{-1}) \geq v(m'm^{-1})$,

(v) for all $j\neq i$, $\{m''\in\mathcal{M}(V)|v(m''m^{-1}) < v(m'm^{-1})\}\cap m'\ZZ[A_{j,a}^{\pm 1}]_{a\in\CC^*}=\emptyset$.

\noindent Then $m'\notin \mathcal{M}(V)$.\end{thm}

To prove this result, we first need some preliminary lemmas.

\subsubsection{Technical lemmas} First let us consider a refined version of the operators $\tau_j$ of \cite{Fre} which allows to study "independently" the subalgebras $\U_q(\Lo\Glie_i)$ of the quantum loop algebra.

 Let $i\in I$, $\Hlie_i^{\perp}=\{\mu\in\Hlie|\alpha_i(\mu)=0\}$ and $A^{(i)}$ be the commutative group 
of monomials generated by variables $Y_{i,a}^{\pm}$ ($a\in\CC^*$), $k_{\mu}$ ($\mu\in\Hlie_i^{\perp}$), 
$Z_{j,c}^{\pm}$ ($j\neq i$, $c\in\CC^*$). Let 
$$\tau_i:A\rightarrow A^{(i)}$$ 
be the group morphism defined by ($j\in I$, 
$a\in\CC^*$): 
$$\tau_i(Y_{j,a})=Y_{j,a}^{\delta_{j,i}}{\prod}_{k\neq 
i,r\in\ZZ}Z_{k,a q^r}^{p_{j,k}(r)}k_{\nu(\Lambda_j)-\delta_{j,i}r_i\alpha_i^{\vee}/2}.$$ 
The $p_{j,k}(r)\in\ZZ$ are defined in the following way : we write 
$\tilde{C}(z)=\frac{\tilde{C}'(z)}{d(z)}$ where $d(z), \tilde{C}'_{j,k}(z)\in\ZZ[z^{\pm}]$ and 
$(D(z)\tilde{C}'(z))_{j,k}={\sum}_{r\in\ZZ}p_{j,k}(r)z^r$. 

\noindent Observe that we have $\nu(\Lambda_j)-\delta_{j,i}r_i\alpha_i^{\vee}/2\in \Hlie_i^{\perp}$ because $\alpha_i(\nu(\Lambda_j)-\delta_{j,i}r_i\alpha_i^{\vee}/2)=\Lambda_j(r_i\alpha_i^{\vee})-r_i\delta_{i,j}=0$.

\noindent This morphism $\tau_i$ was first defined \cite{Fre2}, and then refined in \cite{her06} with the terms $k$ which will be used in the following. Moreover it is proved in \cite[Lemma 3.5]{Fre2} (in \cite[Lemma 20]{her04} with the term $k_0$) that : 

\begin{lem}\label{tauja} For $j\in I$, $a\in\CC^*$, we have $\tau_j(A_{j,a})=Y_{j,aq_j^{-1}}Y_{j,aq_j}k_0$. \end{lem}

\noindent This result indicates that the root monomials $A_{j,a}$ are sent to their analogs of type $sl_2$, as announced above.

\noindent The following result was proved in \cite[Lemma 3.4]{Fre2} without the term $k_{\mu}$, and in \cite[Lemma 21]{her04} the proof was extended for the terms $k_{\mu}$. It gives a decomposition of a $q$-character "compatible" with the action of the subalgebra $\U_q(\Lo\Glie_i)$ :

\begin{lem}\label{aidedeux} Let $V\in\text{Rep}(\U_q(\Lo\Glie))$ and consider a decomposition $\tau_i(\chi_q(V))={\sum}_r P_rQ_r$ where $P_r\in\ZZ[Y_{i,a}^{\pm}]_{a\in\CC^*}$,
$Q_r$ is a monomial in $\ZZ[Z_{j,c}^{\pm}, k_{\lambda}]_{j\neq i,c\in\CC^*, \lambda\in\Hlie_i^{\perp}}$ and all 
monomials $Q_r$ are distinct. Then the
$\U_q(\Lo\Glie_i)$-module $V$ is isomorphic to a direct sum ${\bigoplus}_r V_r$ where
$\chi_q^i(V_r)=P_r$.\end{lem}

The following result gives information on a cyclic $\U_q(\Lo\Glie_j)$-submodule of a $\U_q(\Lo\Glie)$-module : 

\begin{lem}\label{ouun} Let $V\in\text{Rep}(\U_q(\Lo\Glie))$ be a $\U_q(\Lo \Glie)$-module, $m\in\mathcal{M}(L(m))$ and $v\in V_m$. Then for $j\in I$, $\U_q(\Lo\Glie_j).v$ is a sub-$\U_q(\Lo\Hlie)$-module of $V$ and $\chi_q(\U_q(\Lo\Glie_j).v)\in m\ZZ[A_{j,a}^{\pm}]_{a\in\CC^*}$.\end{lem}

\demo From the relation \ref{aidecom}, $\U_q(\Lo\Glie_j).v$ is a sub-$\U_q(\Lo\Hlie)$-module of $V$. Consider the decomposition $\tau_j(\chi_q(V))={\sum}_r P_rQ_r$ of the Lemma 
\ref{aidedeux} and the decomposition of $V$ as a $\U_q(\Lo\Glie_j)$-module:  
$V={\bigoplus}_rV_r$. Then there is $R$ such that $\tau_j(m)$ is a monomial of $P_RQ_R$, and so $v\in V_R$. We have $\U_q(\Lo\Glie_j).v\subset V_R$. Let us write $\tau_j(m)=m_RQ_R$. It follows from \cite[Theorem 7.2]{cm} for $\U_q(\Lo\Glie_j)\simeq \U_{q_j}(\Lo sl_2)$, that the $q$-character of the $\U_q(\Lo\Glie_j)$-module $\U_q(\Lo\Glie_j).v$ is included in $m_R\ZZ[(Y_{j,aq_j^{-1}}Y_{aq_j})^{\pm}]_{a\in\CC^*}$. From Lemma \ref{tauja}, the $q$-character of $\U_q(\Lo\Glie_j).v$ viewed as a $\U_q(\Lo\Hlie)$-module belongs to $m\ZZ[A_{j,a}^{\pm}]_{a\in\CC^*}$.\qed

In the $sl_2$-case, the following Lemma produces a dominant monomial higher than a given monomial in a $q$-character (note that a weak version was proved in \cite[Lemma 3.2 (ii)]{her05}) :

\begin{lem}\label{intsldeux} Let $L$ be a finite dimensional $\U_q(\Lo sl_2)$-module. For $p\in\ZZ$, let $L_p={\sum}_{\{\lambda\in P^*|\lambda(\Lambda^{\vee})\geq p\}}L_{\lambda}$ and $L'_p={\sum}_{r\in\ZZ}x_r^-.L_p$. Then for $m'\in\mathcal{M}(L_p')$ there is $m\in\mathcal{M}(L_p)$ such that 

(i) $m$ is dominant, 

(ii) $m'\leq m$,

(iii) $(\U_q(\Lo sl_2).L_m)\cap L_{m'}\neq \{0\}$.\end{lem}

\demo Let $m'\in\mathcal{M}(L_p')$. Let us prove the result by induction on $\text{dim}(L_p)$ : if $L_p=\{0\}$ we have $L_p'=\{0\}$. In general let $v$ be an $l$-highest weight vector of $L_p$ (it exists, see for example the proof of \cite[Proposition 15]{her04}) and denote by $M$ the corresponding monomial. Consider $V=\U_q(\hat{\Glie}).v$. It is an $l$-highest weight module and so it follows from Theorem \ref{inducb} that $(V_m\neq\{0\}\Rightarrow m\leq M)$. If $V_{m'}\neq \{0\}$ the result is clear with $m=M$. Otherwise consider $L^{(1)}=L/V$. Observe that $\chi_q(L)=\chi_q(V)+\chi_q(L^{(1)})$. We use the induction hypothesis with $L^{(1)}$ and we get $m\in \mathcal{M}((L^{(1)})_p)\subset\mathcal{M}(L_p)$ such that $m\geq m'$ and $(\U_q(\Lo sl_2).(L^{(1)})_m)\cap (L^{(1)})_{m'}\neq \{0\}$. Let $v\in (L^{(1)})_m$ and $\alpha\in \U_q(\Lo sl_2)$ such that $\alpha v\in (L^{(1)})_{m'} - \{0\}$. Let $w\in v + V$ and consider the decomposition $w = w_m + w'$ where $w_m\in L_m$ and $w'\in \bigoplus_{m''\neq m} L_{m''}$. Consider $v\in (L^{(1)})_m$, we have $w'\in V$ and $w_m\in v + V$. Then $\alpha w_m = v' + v''\in L_{m'}\oplus V$ where $v'\neq 0$. As $V_{m}' = \{0\}$, there is $h\in\U_q(\Hlie)$ such that $h\alpha w_m=hv'\neq 0$ and so we get the result.\qed

An analog result is available for general type :

\begin{lem}\label{oudeux} Let $V = L(m)$ be a $\U_q(\Lo \Glie)$-module simple module and $m'<m$ in $\mathcal{M}(L(m))$. Then there is $j\in I$ and $M'\in\mathcal{M}(V)$ such that

(i) $M'$ is $j$-dominant,

(ii) $M'>m'$, 

(iii) $M'\in m'\ZZ[A_{j,b}]_{b\in\CC^*}$,

(iv) $((\U_q(\Lo\Glie_j).V_{M'})\cap (V)_{m'})\neq \{0\}$.\end{lem}

(A weak version of the following lemma was proved in the proof of \cite[Lemma 4.4]{her06} with different notations).

To prove this result, we need the following additional notations : for $M\in A^{(i)}$, we define $\mu(M)\in\Hlie_i^{\perp}$, $u_{i,a}(M)\in\ZZ$, by : 
$$M\in k_{\mu(M)}{\prod}_{a\in\CC^*}Y_{i,a}^{u_{i,a}(M)}\ZZ[Z_{j,c}^{\pm}]_{j\neq i, c\in\CC^*}.$$
We also set $u_i(M)={\sum}_{a\in\CC^*}u_{i,a}(M)$. Observe that for $m\in A$ and $a\in\CC^*$ we have $u_{i,a}(m)=u_{i,a}(\tau_i(m))$ and :
$$\nu(\omega(m))=\mu(\tau_i(m))+u_i(m)r_i\alpha_i^{\vee}/2=\mu(\tau_i(m))+u_i(\tau_i(m))r_i\alpha_i^{\vee}/2,$$
or equivalently 
$$\mu(\tau_i(m))=\nu(\omega(m))-\alpha_i(\nu(\omega(m)))\alpha_i^{\vee}/2.$$
(See the definition of \cite[Section 5.5]{her04}.) Now let us prove Lemma \ref{oudeux} :

\demo For $m''\in\mathcal{M}(V)$ denote $w(m'')=v(\omega(m'')-\omega(m))$. Let 
$$W=\bigoplus_{\{m''|w(m'')<w(m')\}}V_{m''}.$$ 
As $V$ is an $l$-highest weight module, there is $j\in I$ such that $(\U_q(\Lo\Glie_j).W)_{m'}\neq \{0\}$. Consider the decomposition $\tau_j(\chi_q(V))={\sum}_rP_rQ_r$ of Lemma \ref{aidedeux} and the decomposition of $V$ as a $\U_q(\Lo\Glie_j)$-module:  $V = {\bigoplus}_rV_r$. 

For a given $r$, consider $M_r\in\mathcal{M}(V)$ such that 
$\tau_j(M_r)$ appears in $P_rQ_r$. For another such $M$, we have $\mu(\tau_j(M))=\mu(\tau_j(M_r))$ and so 
$$\omega(MM_r^{-1})=u_j(\tau_j(MM_r^{-1}))\alpha_j^{\vee}/2,$$ 
and 
\begin{equation*}
u_j(\tau_j(M))=u_j(\tau_j(M_r))-2w(M)+2w(M_r)=2(p-w(M))+p_r,
\end{equation*} 
where $p_r=-2p+2w(M_r)+u_j(\tau_j(M_r))$ (it does not 
depend of $M$). So we have $w(M)\leq p\Leftrightarrow u_j(\tau_j(M))\geq p_r$. So 
$W={\bigoplus}_r((V_r)_{\geq p_r})={\bigoplus}_r(V_r\cap W)$. As $V_r$ is a sub 
$\U_q(\Lo\Glie_j)$-modules of $V$, we have $W_j={\bigoplus}_r(V_r\cap W_j)$. 
Let $M\in\mathcal{M}(W_j)$ and $R$ such that $\tau_j(M)$ is a monomial of $P_RQ_R$. We can apply Lemma \ref{intsldeux} to the $\U_{q}(\Lo\Glie_j)$-module $V_R$ with $p = p_R$ and the monomial $Q_R^{-1}\tau_j(m)$ : we get $m''\in\mathcal{M}(V_R)$ dominant such that $Q_R^{-1}\tau_j(M)\in m'' \ZZ[(Y_{j,a}Y_{j,aq_j^2})^{-1}]_{a\in\CC^*}$ and $((\U_q(\Lo\Glie_j).(V_R)_{m''})\cap (V_R)_{Q_R^{-1}\tau_j(M)})\neq \{0\}$. Let us translate this result in terms of monomials of $\chi_q(V)$. Consider the $j$-dominant monomial $M'=\tau_j^{-1}(Q_Rm')$. Then $M'\in\mathcal{M}(W)$ and $(\U_q(\Lo\Glie_j).V_{M'})\cap V_M\neq \{0\}$. From Lemma \ref{tauja} we have $M\in M'\ZZ[A_{j,b}^{-1}]_{b\in\CC^*}$.\qed

\subsubsection{Proof of Theorem \ref{racourc}}

Suppose that $m'\in \mathcal{M}(V)$. Let 
$$W=\bigoplus_{\{M'\leq m|v(M'm^{-1}) < v(m'm^{-1})\}}V_{M'}.$$ 
As $V$ is an $l$-highest weight modules, there is $k\in I$ such that $\sum_{r\in\ZZ}(x_{k,r}^-.W)_{m'}\neq \{0\}$. From condition (v) and Lemma \ref{ouun}, we have $k=i$. From Lemma \ref{oudeux} and condition (i), we have $(V_{m'}\cap\U_q(\Lo\Glie_i) V_M)\neq\{0\}$. Consider $u\in V_M$ and $x\in (\U_q(\Lo\Glie_i).u\cap V_{m'})$ such that $x\neq 0$. From condition (ii), $u$ is an highest weight vectors for $\U_q(\Lo\Glie_i)$, so $x\in \sum_{r\in\ZZ}\CC x_{i,r}^-.u$. By condition (iii), $x$ is in the maximal proper $\U_q(\Lo\Glie_i)$-submodule of $\U_q(\Lo\Glie_i).x$. By condition $(iv)$, $v(m'M^{-1})$ is maximal for this condition. So for all $r\in\ZZ$, we have $x_{i,r}^+(x)=0$. For $j\neq i$, $r\in\ZZ$, it follows from Lemma \ref{ouun} that $x_{j,r}^+(x)\in \bigoplus_{m''\in m'\ZZ[A_{j,a}^{\pm 1}]_{a\in\CC^*}}V_{m''}$, and so from condition (v) we have $x_{j,r}^+(x)=0$. So $\U_q(\Lo\Glie).x$ is a proper submodule of $V$, contradiction.\qed

\subsection{Other preliminary results}\label{secoth}

In this section, $\Glie$ is an arbitrary semi-simple Lie algebra. We prove additional preliminary results.

\subsubsection{$q$-characters of simple modules} 

\begin{lem}\label{produit} Let $L(m_1), L(m_2)$ be two simple modules. Then $L(m_1m_2)$ is a subquotient of $L(m_1)\otimes L(m_2)$. In particular $\mathcal{M}(L(m_1m_2))\subset\mathcal{M}(L(m_1))\mathcal{M}(L(m_2))$.\end{lem}

This first part of the lemma is proved in \cite{Cha2}, and the second part is direct from \cite{Cha2, Fre}. 

As a direct consequence of Theorem \ref{inducb}, we have :

\begin{lem}\label{plusgrand} Let $a\in\CC^*$ and $m$ be a monomial of $\ZZ[Y_{i,aq^r}]_{i\in I,r\geq 0}$. Then for $m'\in\mathcal{M}(L(m))$ and $b\in\CC^*$, ($v_{i,b}(m'm^{-1})\neq 0\Rightarrow b\in aq^{r_i + \NN}$).\end{lem}

\noindent (Observe that it also a direct consequence of Lemma \ref{fundless} as a simple module is a subquotient of a tensor product of fundamental representations.)

The following result gives information on the sub $\U_q(\Lo\Glie^J)$-module generated by an highest weight vector (the definition of $L_J(m)$ and $L^J(m^{\rightarrow(J)})$ were given in section \ref{compred}) :

\begin{lem}\label{depart} Let $m$ be a dominant monomial and $J\subset I$. Let $v$ be an highest weight vector of $L(m)$ and $L'\subset L(m)$ the $\U_q(\Lo \Glie_J)$-submodule of $L(m)$ generated by $v$. Then $L'$ is an $\U_q(\Lo\Hlie)$-submodule of $L(m)$ and $\chi_q(L')=L_J(m)$.\end{lem}

In particular for $\mu\in \omega(m)-{\sum}_{j\in J}\NN\alpha_j$, we have $$\text{dim}((L(m))_\mu)=\text{dim}((L^J(m^{\rightarrow(J)}))_{\mu^{\rightarrow(J)}}),$$ 
where $\mu^{\rightarrow(J)}={\sum}_{j\in J}\mu(\alpha_j^{\vee})\omega_j$.

\demo  It is clear that $L'=\bigoplus_{\mu\in \omega(m)-\sum_{j\in J}\NN.\alpha_j}(L(m))_\mu$. So it is an $\U_q(\Lo \Hlie)$-submodule of $L(m)$ and $\chi_q(L')$ makes sense. Moreover $\chi_q(L')\in m\ZZ[A_{j,a}^{\pm}]_{j\in J,a\in\CC^*}$ and $m'\in\mathcal{M}(L')$ is uniquely determined by $(m')^{\rightarrow(J)}$. So it suffices to prove that $L'\simeq L^J(m^{\rightarrow(J)})$ as $\U_q(\Lo\Glie_J)$-module. As $L'$ is an highest weight $\U_q(\Lo \Glie_J)$-module of highest weight monomial $m^{\rightarrow(J)}$, it suffices to prove that $L'$ is simple. If it is not simple, there is $w\in L'\cap (L(m))_\mu$ where $\mu < \omega(m)$ and such that for all $j\in J$, $m\in\ZZ$, $x_{j,m}^+.w = 0$. But as $L(m)$ is an highest weight module and the weight of $L'$ are in $\omega(m) - \sum_{j\in J} \NN \alpha_j$, for weight reason we have :
$$\forall j\in (I-J)\text{ , }\forall m\in\ZZ\text{ , }x_{j,m}^+(L')=\{0\}.$$ 
So $\U_q(\Lo \Glie).w$ is a proper submodule of $L(m)$, contradiction.\qed

\subsubsection{Thin modules and thin monomials}

Let us introduce the notion of thin module :

\begin{defi} A $\U_q(\Lo\Glie)$-module $V$ is said to be thin if his $l$-weight spaces are of dimension $1$.\end{defi}

In \cite[Theorem 3.2]{her05}, we proved that for $\Glie$ of type $A$, $B$, $C$, all fundamental representations are thin (this result was also proved later by a different method in \cite{cm2}. It should be also possible to check this result directly from the formulas in \cite{ks}). We will discuss in more detail the notion of thin modules in \cite{miniaff}, but let us give some results that will be used in the present paper.

\begin{lem}\label{thinmon} Let $V$ be a $\U_q(\Lo\Glie)$-module and $m'\in\mathcal{M}(V)$ such that there is $i\in I$ satisfying $\text{Min}\{u_{i,a}(m')|a\in\CC^*\}\leq -2$. Then there is $M\in\mathcal{M}(V)$ such that

$M > m'$,

$M$ is $i$-dominant,

$\text{Max}\{u_{i,b}(M)|b\in\CC^*\} \geq 2$.\end{lem}

\demo Consider $L_i(M)$ occurring in the decomposition of $\chi_q(V)$ described in Proposition \ref{jdecomp} and such that $m'$ is a monomial of $L_i(M)$. $L_i(M)$ corresponds to the $q$-character $\chi_q^i(W)$ where $W$ is a $\U_q(\Lo\Glie_i)$-simple module, so subquotient of a standard module. In particular $m'$ appears in
\begin{equation}\label{formu}M \prod_{a\in\CC^*} (Y_{i,a}^{-1}(1 + A_{i,aq_i}^{-1}))^{u_{i,a}(M)}.\end{equation} 
By hypothesis there is $b\in\CC^*$ such that $u_{i,b}(m') \leq -2$. As $m'$ appears in the formula (\ref{formu}), necessarily $(1+A_{i,aq_i}^{-1})$ appears at least twice in (\ref{formu}), and so $u_{i,bq_i^{-1}}(M)\geq 2$. Moreover by construction $M > m'$ and $M$ is $i$-dominant.\qed

\begin{defi} A monomial $m$ is said to be thin if $\text{Max}_{i\in I, a\in\CC^*} |u_{i,a}(m)|\leq 1$.\end{defi}

\begin{lem}\label{dimun} Let $V$ be a special module such that 
$$\text{Max}\{u_{i,a}(m)|m\in\mathcal{M}(V), i\in I, a\in\CC^*\}\leq 1.$$ 
Then $V$ is thin. Moreover all $m\in\mathcal{M}(V)$ are thin. \end{lem}

\demo In \cite[Proposition 3.3]{her05}, the first statement is proved for fundamental representations.
The proof of the first statement of the Lemma is the same ($\chi_q(V)$ is given by the Frenkel-Mukhin algorithm, and so the property is proved by induction on the weight of monomials, see the proof of \cite[Proposition 3.3]{her05} for details). 

\noindent Now consider $m'\in\mathcal{M}(V)$. If $m'$ is not thin, there is $i\in I$ and $a\in\CC^*$
such that $u_{i,a}(m')\leq -2$. From Lemma \ref{thinmon}, there is another monomial 
$M\in\mathcal{M}(V)$ and $b\in\CC^*$ such that $u_{i,b}(M)\geq 2$, contradiction with
the hypothesis on $V$, so $m'$ is thin.\qed

\begin{prop}\label{proofthin} If $V$ is thin then all $m\in\mathcal{M}(V)$ are thin. If $V$ is special and all $m\in\mathcal{M}(V)$ are thin, then $V$ is thin.\end{prop}

\demo If $V$ is special and all $m\in\mathcal{M}(V)$ are thin, then the hypothesis of Lemma \ref{dimun} are satisfied and so $V$ is thin. 

\noindent For the first statement, suppose that $V$ is thin and that there is a monomial of $\mathcal{M}(V)$ which is not thin. We can suppose there is $m\in\mathcal{M}(V)$, $i\in I$, $a\in\CC^*$ such that $u_{i,a}(m)\geq 2$ (in the case $u_{i,a}(m)\leq -2$ it follows from Lemma \ref{thinmon} that there is another monomial satisfying the condition with $\geq 2$). Consider $L_i(M)$ occurring in the decomposition of $\chi_q(V)$ described in Proposition \ref{jdecomp} and such that $m$ is a monomial of $L_i(M)$. We can see as in the proof of Lemma \ref{thinmon} that there is $b\in\CC^*$ satisfying $u_{i,b}(M)\geq 2$. From the explicit description of simple modules in Proposition \ref{aidesldeux} in the case $sl_2$, the monomial $MA_{i,bq_i}^{-1}\prod_{r > 0}A_{i,b(q_i)^{r+1}}^{-u_{i,bq_i^r}(M)}$ occurs with multiplicity at least $2$ in the $q$-character of the $\U_{q_i}(\Lo sl_2)$-module $L(M^{\rightarrow(i)})$, and so it is not a thin module. As the coefficients in the decomposition of Proposition \ref{jdecomp} are positive, there is an $l$-weight space of $V$ of dimension at least $2$, and so $V$ is not thin.\qed

\begin{lem}\label{etoile} Let $L(m)$ be a simple $\U_q(\Lo\Glie)$-module and $(m',i)\in\mathcal{M}(L(m))\times I$ such that 

all $m''\in\mathcal{M}(L(m))$ satisfying $v(m''m^{-1}) < v(m'm^{-1})$ is thin, 

$m'$ is not $i$-dominant. 

\noindent Then there is $a\in\CC^*$ such that $u_{i,a}(m') < 0$ and $m'A_{i,aq_i^{-1}}\in\mathcal{M}(L(m))$.\end{lem}

\demo Consider $L_i(M)$ occurring in the decomposition of $\chi_q(V)$ described in Proposition \ref{jdecomp} and such that $m$ is a monomial of $L_i(M)$. From the first hypothesis $M$ is thin. If $L_i(M)$ correspond to a Kirillov-Reshetikhin module of type $sl_2$, the result follows from the explicit formula of Proposition \ref{aidesldeux} (1). In general $L_i(M)$ is also known from the explicit description of $q$-characters of simple modules in the $sl_2$-case in Proposition \ref{aidesldeux} (3), and $L_i(M)$ corresponds to a product of Kirillov-Reshetikhin modules $L_i(M) = \prod_k W_k$. As $M$ is thin we have moreover the following property : for $m_1$ appearing in $W_k$ and $m_2$ appearing in $W_{k'}$, we have
$$u_{i,a}(m_1) \neq 0 \text{ and }u_{i,a}(m_2)\neq 0\Rightarrow k = k'.$$
And so the result can be reduced to the case of Kirillov-Reshetikhin modules.\qed

\begin{lem}\label{stara} Suppose that $\Glie = sl_{n+1}$ and $L(m)$ be a simple $\U_q(\Lo\Glie)$-module. 
Let $(m',i,a )\in\mathcal{M}(L(m))\times I\times \CC^*$ such that :

all $m''\in\mathcal{M}(L(m))$ satisfying $v(m''m^{-1}) < v(m'm^{-1})$ is thin, 

$u_{i,a}(m') = -1$,

$m'Y_{i,a}$ is dominant. 

\noindent Then there is $M\in\mathcal{M}(L(m))$ dominant such that $M > m'$ and $v_n(m'M^{-1})\leq 1$, $v_1(m'M^{-1})\leq 1$.\end{lem}

\demo By using Lemma \ref{etoile}, we construct inductively a sequence of monomials of $\mathcal{M}(L(m))$ starting with $m'$. Indeed as $u_{i,a}(m') = -1$ we first get $m'A_{i,aq^{-1}}\in \mathcal{M}(L(m))$. Then from the property $m'Y_{i,a}$ dominant we have 
$$(u_{i-1,b}(m'A_{i,aq^{-1}}) < 0\text{ or }u_{i+1,b}(m'A_{i,aq^{-1}}) < 0) \Rightarrow b = aq^{-1}.$$ 
Then we use again Lemma \ref{etoile} $(i-1,aq^{-1})$ and $(i+1,aq^{-1})$ when it is possible. We get a monomial and we apply Lemma \ref{etoile} with $(i-2,aq^{-2})$ and $(i+2,aq^{-2})$ when it is possible. We continue by induction until this is not possible, and we get a monomial : 
\begin{equation*}
\begin{split}
m_1 = &m'(A_{i,aq^{-1}}A_{i-1,aq^{-2}}\cdots A_{i-\alpha,aq^{-1-\alpha}})
\\&\times (A_{i+1,aq^{-2}}A_{i+2,aq^{-3}}\cdots A_{i+\beta,aq^{-1-\beta}})\in\mathcal{M}(L(m)),
\end{split}
\end{equation*}
where $\alpha,\beta\geq 0$, $i-\alpha\geq 1$, $i+\beta\leq n$. By construction $m_1$ is $(I-\{i\})$-dominant and we have $(u_{i,b}(m_1) <0\Rightarrow b=aq^{-2})$. If $\alpha = 0$ or $\beta =0$, $m_1$ is dominant and we take $M=m_1$. Otherwise, we can suppose $\alpha\geq \beta$ (the case $\beta \geq\alpha$ can be treated in the same way). As at each step we get by construction thin monomials, we continue by induction, and for $2\leq r\leq \beta +1$, we have
\begin{equation*}
\begin{split}
m_r=&m_{r-1}(A_{i,aq^{1-2r}}A_{i-1,aq^{-2r}}\cdots A_{i-\alpha +r-1,aq^{-\alpha-r}})
\\&\times (A_{i+1,aq^{-2r}}A_{i+2,aq^{-2r-1}}\cdots A_{i+\beta -r+1,aq^{-r-\beta}})\in\mathcal{M}(L(m)),
\end{split}
\end{equation*}
and $m_r$ is $(I-\{i\})$ dominant. Moreover $m_{\beta +1}$ is dominant, so we take $M=m_{\beta +1}$. By construction we have $M > m'$ and $v_n(m'M^{-1})\leq 1$, $v_1(m'M^{-1})\leq 1$.\qed

\begin{lem}\label{remont} Let $\Glie = sl_{n+1}$ and $L(m)$ be a simple $\U_q(\Lo\Glie)$-module. Let $(m',j)\in\mathcal{M}(L(m))\times I$ such that 

all $m''\in\mathcal{M}(L(m))$ satisfying $v(m''m^{-1}) < v(m'm^{-1})$ is thin,

$m'$ is $(I-\{j\})$-dominant,

if $j\leq n-1$, then for all $a\in\CC^*$, $(u_{j,a}(m')<0\Rightarrow u_{j+1,aq^{-1}}(m') > 0)$. 

\noindent Then there is $M\in\mathcal{M}(L(m))$ dominant of the form 
$$M = m'\prod_{\{a\in\CC^*|u_{j,a}(m')<0\}}(A_{j,aq^{-1}}A_{j-1,aq^{-3}}\cdots A_{i_a,aq^{i_a-j-1}}),$$
where for $a\in\CC^*$, $1\leq i_a\leq j$.
\end{lem}

\demo If $j < n$, the additional hypothesis $(u_{j,a}(m')<0\Rightarrow u_{j+1,aq^{-1}}(m') > 0)$ allows to use the result for $\Glie_{\{1,\cdots,j\}}$. So we can suppose that $j = n$. We prove the result by induction on $n$. For $n=1$ the result is clear. In general, by using Proposition \ref{jdecomp} we get $m_1\in\mathcal{M}(L(m))$ $n$-dominant such that $m'$ is a monomial of $L_n(m_1)$. As $m_1 > m'$ and $v(m_1m^{-1})\leq P$, we have by the explicit description of $L_n(m_1)$ in Proposition \ref{aidesldeux} : 
$$m'(m_1)^{-1}=\prod_{\{a\in\CC^*|u_{n,a}(m')<0\}}A_{n,aq^{-1}}.$$
Moreover by construction : 

$m_1$ is $\{1,\cdots,n-2\}$-dominant,

$\forall a\in\CC^*, (u_{n-1,a}(m_1) < 0 \Rightarrow (u_{n,aq^{-1}}(m_1) = 1\text{ and }u_{n-1,a}(m_1) = -1))$.

\noindent By Lemma \ref{etoile} there is $m_2\in\mathcal{M}(L(m))$ which is $\{1,\cdots,n-1\}$-dominant and such that $m_1$ is a monomial of $L_{\{1,\cdots, n-1\}}(m_2)$. Then by using the induction property for $\Glie_{\{1,\cdots,n-1\}}$ on $m_1$ monomial of $L_{\{1,\cdots,n-1\}}(m_2)$, we get the monomial $M$.\qed

\section{Proof of Theorem \ref{mainres}}\label{small}

In this section, $\Glie$ is simply-laced. We complete the proof of Theorem \ref{mainres} : after a technical lemma on dominant monomials (Lemma \ref{lequeldom}), fundamental representations (Proposition \ref{fundsmall}) and standard modules of the form $M(X_{2,a}^{(i)})$ (Proposition \ref{kegaldeux}) are studied. Then the type $A$ is discussed (Proposition \ref{equismalla}), and finally we give the proof of Theorem \ref{mainres} for the general case.

\subsection{Dominant monomials} First let us prove some properties of dominant monomials lower than a monomial $X_{k,a}^{(i)}$. To do this, let us define the following number attached to the structure of the Dynkin diagram : for $i,j\in I$, we denote by $d(i,j)$ the minimal $d$ such that there is a sequence $(i_1,\cdots,i_d)\in I^d$ satisfying $i=i_1$, $j=i_d$ and for all $k\in\{1,\cdots,d-1\}$, $C_{i_k,i_{k+1}}=-1$.

\begin{lem}\label{lequeldom} Let $i\in I$, $a\in\CC^*$, $k\geq 0$ and $m = X_{k,a}^{(i)}$. Let $m'\leq m$ dominant. Then we have :

$m'm^{-1}\in\ZZ[A_{j,aq^l}^{-1}]_{j\in I,l\in\ZZ}$,

$\forall j\in I, l\in\ZZ$, $v_{j,a q^l}(m'm^{-1}) > 0\Rightarrow (d(i,j)+1-k\leq l\leq k-1-d(i,j))$,
 
$\forall j\in I$, $v_j(m'm^{-1}) > 0\Rightarrow d(i,j)\leq k-1$.\end{lem}

\demo The last statement in a direct consequence of the second statement. 

\noindent Let us prove that for any $j\in I$, $b\in\CC^*$ we have :
$$v_{j,b}(m'm^{-1})\neq 0 \Rightarrow b\in aq^{k-1-d(i,j)-\NN}.$$ 
We prove this statement by induction on $d(i,j)$. 

\noindent For $d(i,j) = 1$, we have $j = i$. Suppose that there is $b\in(\CC^* - aq^{k-2-\NN})$ such that $v_{i,b}(m'm^{-1}) > 0$. Let $L\in\ZZ$ maximal such that there is $p\in I$ satisfying $v_{p,b q^L}(m'm^{-1}) > 0$. We have $bq^L\notin aq^{k-2-\NN}$. As $m' \neq m$, we have $m' < m$ and $m'm^{-1}$ is right negative. So $u_{p,bq^{L+1}}(m'm^{-1}) < 0$. As moreover $u_{p,c}(m) > 0$ implies $c\in aq^{k-1-\NN}$, we have $u_{p,bq^{L+1}}(m') = u_{p,bq^{L+1}}(m'm^{-1}) < 0$. So $m'$ is not dominant, contradiction. 

\noindent In general suppose that $d(i,j) \geq 2$ and that there is $b\in(\CC^* - aq^{k-d(i,j)-1 -\NN})$ such that $v_{j,b}(m'm^{-1}) > 0$. If $b\notin aq^\ZZ$, we can prove as in the previous case that $m'$ is not dominant, contradiction. Otherwise let $L$ maximal such that 
$$\sum_{\{p\in I|d(i,p)\geq d(i,j)\}}v_{p,aq^L}(m'm^{-1}) > 0.$$ 
As $v_{j,b}(m'm^{-1}) > 0$, we have $L > k-1-d(i,j)$. Let $P\in I$ such that $d(i,P)\geq d(i,j)$ and $v_{P,aq^L}(m'm^{-1}) > 0$. We have 
$$u_{P,aq^{L+1}}(\prod_{\{p\in I|d(i,p)\geq d(i,j)\}}\prod_{c\in\CC^*}A_{p,c}^{-1})< 0.$$
As $u_{P,aq^{L+1}}(m) = 0$ and $m'$ is dominant, there is $j'$ satisfying $d(i,j') = d(i,j) - 1$ and $v_{j',aq^{L+1}}(m'm^{-1}) > 0$. But $L+1 \geq k  +1 - d(i,j) = k - d(i,j')$, contradiction with the induction hypothesis. 

\noindent In the same way we can prove that for any $j\in I$, $b\in\CC^*$ : 
$$v_{j,b}(m'm^{-1})\neq 0 \Rightarrow b\in aq^{-k+1+d(i,j) + \NN}.$$ 
This implies the first two statements of the Lemma.\qed

\subsection{Fundamental representations and $k = 2$ case}

\begin{prop}\label{fundsmall} All fundamental representations are small.\end{prop}

\demo Let $i\in I$ and $a\in\CC^*$. Then from Lemma \ref{lequeldom}, a monomial satisfying $m'<Y_{i,a}$ is not dominant. So $V_i(a)$ is small.\qed

\begin{prop}\label{kegaldeux} Let $i\in I$, $a\in\CC^*$. Then $M(X_{2,a}^{(i)})$ is small.\end{prop}

\demo 
From Lemma \ref{lequeldom}, a dominant monomial $m' < X_{2,a}^{(i)}$ is equal to 
$$m' = X_{2,a}^{(i)}A_{i,a}^{-1} = \prod_{j\in I|C_{i,j} = -1} Y_{j,a}.$$ 
Consider a monomial $m'' < m'$. 

\noindent Suppose that there is $j\in I$, $b\in(\CC^*-aq^\ZZ)$ such that $v_{j,b}(m''(m')^{-1}) > 0$. Let $L\in\ZZ$ maximal such that there is $p\in I$ satisfying $v_{p,b q^L}(m''(m')^{-1}) > 0$. We have $u_{p,bq^{L+1}}(m'') = u_{p,bq^{L+1}}(m''(m')^{-1}) < 0 $ and so $m''$ is not dominant. 

\noindent Otherwise let $L\in\ZZ$ maximal such that there is $p\in I$ satisfying $v_{p,a q^L}(m''(m')^{-1}) > 0$. If $L\geq 0$, we can prove as in the previous case that $m''$ is not dominant. Otherwise let $L' < 0$ minimal such that there is $p\in I$ satisfying $v_{p,a q^{L'}}(m''(m')^{-1}) > 0$. We have $u_{p,bq^{L' - 1}}(m'') = u_{p,bq^{L'-1}}(m''(m')^{-1}) < 0 $, and so $m''$ is not dominant.

\noindent So $L(m')$ is special and $M(X_{2,a}^{(i)})$ is small.\qed

\subsection{Type $A$} In this section $\Glie$ is of type $A$. 

\begin{prop}\label{equismalla} Let $k\geq 1, i\in I, a\in\CC^*$. Then $M(X_{k,a}^{(i)})$ is small if and only if ($i = 1$ or $i = n$ or $k\leq 2$).\end{prop} 

\noindent In particular for $\Glie = sl_2$ or $\Glie = sl_3$, all $M(X_{k,a}^{(i)})$ are small.

\noindent We prove this proposition in three steps : 

(1) we determine the dominant monomials $m'$ such that $m'\leq X_{k,a}^{(1)}$ (Lemma \ref{formun}), 

(2) we prove that the corresponding simple modules are special (Proposition \ref{spesmalllem}),

(3) we study the remaining cases (Lemma \ref{reste}). 

\begin{lem}\label{formun} Let $k\geq 1, a\in\CC^*$ and $m'\leq X_{k,a}^{(1)}$ dominant. Then $m'$ is of the form 
$$m'=Y_{i_1,aq^{l_1}}Y_{i_2,aq^{l_2}}\cdots Y_{i_R,aq^{l_R}},$$
where $R\geq 0$, $i_1 ,i_2,\cdots, i_R\in I$, $l_1,l_2,\cdots , l_R\in \ZZ$ satisfy for all $1\leq r\leq R - 1$ :
$$l_{r+1}-l_r \geq i_r+i_{r+1}.$$
\end{lem}

\demo Let $m = Y_{1,a}Y_{1,aq^2}\cdots Y_{1,aq^{2(k-1)}}$ and $m'\leq m$ dominant. For $i\in I,l\in\ZZ$, let us denote $v_{i,l}=v_{i,aq^l}(m'm^{-1})$ and $u_{i,l}=u_{i,aq^l}(m')$. We denote $v_{n+1,l}=0$. As $m'$ is dominant, we have for $2\leq i\leq n$ and $l\in\ZZ$ :
$$v_{i,l-1}+v_{i,l+1}\leq v_{i-1,l}+v_{i+1,l},$$
$$v_{1,l-1}+v_{1,l+1}\leq 1 + v_{2,l}.$$
From Lemma \ref{lequeldom}, for $l\leq i - 1$ or $l\geq 2k - i - 2$, we have $v_{i,l}=0$. 
Let us prove that for all $i\in I$, $(v_{i,l}\neq 0\Rightarrow (l\in i + 2\ZZ))$. Indeed $m''=\prod_{i\in I, l\in i+1+2\ZZ}A_{i,aq^l}^{-v_{i,l}}$ is right negative and for all $i\in I$, $l\in i+2\ZZ$, $u_{i,l}(m'')=u_{i,l}(m')$. So $m''=1$.

Let us prove that for all $l\in \ZZ$ we have $v_{1,l}\leq 1$, and for all $n\geq i\geq 2$, $l\in\ZZ$ we have $v_{i,l}\leq v_{i-1,l-1}$. We prove the result by induction on $d=l-i\geq 0$. First suppose that $d=0$. First we have $v_{1,1}\leq -v_{1,-1}+1+v_{2,1}=1$. For $i\geq 2$, $v_{i,i}\leq -v_{i,i-2}+v_{i-1,i-1}+v_{i+1,i-1}=v_{i-1,i-1}$. Now consider a general $d > 0$. First we have $v_{1,1+d}\leq 1+v_{2,d}-v_{1,d-1}$. But by the induction hypothesis, $v_{2,d}\leq v_{1,d-1}$. So $v_{1,1+d}\leq 1$. For $i\geq 2$, $v_{i,i+d}\leq (v_{i+1,d+i-1}-v_{i,d+i-2})+v_{i-1,d+i-1}$. But by the induction hypothesis, $v_{i+1,d+i-1}-v_{i,d+i-2}\leq 0$, and so $v_{i,i+d}\leq v_{i-1,d+i-1}$.

In particular for all $i\in I, l\in\ZZ$, $v_{i,l}\leq 1$.

In the same way, for all $n\geq i\geq 2$, $v_{i,l}\leq v_{i-1,l+1}$. Let $n\geq i\geq 2$. We have proved $v_{i,l}\leq \text{Min}\{v_{i-1,l-1},v_{i-1,l+1},1\}$. In particular 
$$(v_{i,l}=1\Rightarrow v_{i-1,l-1}=v_{i-1,l+1}=1).$$ 
Moreover if $v_{i,l-1}=v_{i,l+1}=1$, we have $2=v_{i,l-1}+v_{i,l+1}\leq v_{i+1,l}+v_{i-1,l}$ and so $v_{i+1,l}=v_{i-1,l}=1$. So 
$$(v_{i,l}=1\Leftrightarrow v_{i-1,l-1}=v_{i-1,l+1}=1).$$ 
As a conclusion, this can be rewritten in the following way. $m'm^{-1}$ is of the form : 
$$m'm^{-1}=B_{p_1,f_1}B_{p_2,f_2}\cdots B_{p_R,f_R},$$ 
where $R\geq 0$, $n-1\geq p_1,\cdots, p_R\geq 0$, $f_1,\cdots, f_R\in\ZZ$, 
\begin{equation*}
\begin{split}
B_{p,f}=&(A_{1,aq^{f-p}}A_{1,aq^{f+2-p}}\cdots A_{1,aq^{f+p}})
\\&\times(A_{2,aq^{f+1-p}}A_{2,aq^{f+3-p}}\cdots A_{2,aq^{f+p-1}})\cdots (A_{p+1,aq^f}),
\end{split}
\end{equation*}
$f_i-p_i\in 1+2\ZZ$, $f_1-p_1\geq 1$, $f_R+p_R\leq 2k-3$ and $f_i+p_i + 4 \leq f_{i+1}-p_{i+1}$.

If $p\leq n-2$, we have 
$$B_{p,l}=(Y_{1,q^{f-p-1}}^{-1}Y_{1,q^{f-p+1}}^{-1}\cdots Y_{1,q^{f+p+1}}^{-1})Y_{p+2,aq^f},$$ 
and we have
$$B_{n-1,l}= Y_{1,q^{f-n}}^{-1}Y_{1,q^{f-n+2}}^{-1}\cdots Y_{1,q^{f+n}}^{-1}.$$
So we get the result.\qed

\begin{prop}\label{spesmalllem} Let $m = Y_{i_1,aq^{l_1}}Y_{i_2,aq^{l_2}}\cdots Y_{i_R,aq^{l_R}}$ where $R\geq 0$, $i_1,i_2,\cdots, i_R\in I$, $l_1,l_2,\cdots , l_R\in \ZZ$ satisfying for all $1\leq r\leq R-1$, $l_{r+1}-l_r \geq i_r+i_{r+1}$. Then :

(1) For $m'\in \mathcal{M}(L(m))$, if $v_{i_R,aq^{l_R-1}}(m')\geq 1$ then $v_{i_R,aq^{l_R+1}}(m')\geq 1$.

(2) $L(m)$ is special.

(3) $L(m)$ is thin.
\end{prop}

To prove this Proposition, we will need the following direct consequence of the results in \cite{Fre2} :

\begin{lem}\label{fundalgo} Let $V$ be a fundamental representation of a quantum loop algebra $\U_q(\Lo\Glie)$ and let $Y_{i,a}$ (resp. $Y_{j,b}^{-1}$) be the highest (resp. lowest) weight monomial of $\chi_q(V)$. Then we have :
$$\chi_q(V)\in Y_{i,a} (1 + A_{i,aq}^{-1} (1 + \sum_{\{k\in I|C_{i,k} = -1\}} A_{k,bq^2}^{-1} . \ZZ[A_{l,d}^{-1}]_{l\in I,d\in\CC^*})),$$
$$\chi_q(V)\in Y_{j,b}^{-1} ( 1 + A_{j,bq^{-1}} (1 + \sum_{\{k\in I|C_{j,k} = -1\}} A_{k,bq^{-2}} . \ZZ[A_{l,d}]_{l\in I,d\in\CC^*})).$$
\end{lem}

\demo As $V$ is special, we can use the algorithm proposed by Frenkel-Mukhin \cite{Fre2} to compute $\chi_q(V)$ (see \cite[Section 5.5]{Fre2} for details) : we start with $Y_{i,a}$. Then we get $Y_{i,a}A_{i,aq}^{-1}$ with multiplicity $1$ as $L_i(Y_{i,a}) = Y_{i,a} + Y_{i,a}A_{i,aq}^{-1}$. As 
$$Y_{i,a}A_{i,aq}^{-1} = Y_{i,aq^2}^{-1}\prod_{\{k\in I|C_{i,k} = -1\}}Y_{k,aq},$$ 
the next step of the algorithm gives the monomials $Y_{i,a}A_{i,aq}^{-1}A_{k,aq^2}^{-1}$ with multiplicity one, and then inductively the other monomials occurring in $\chi_q(V)$ are lower then these monomials.

The second statement is obtained by the duality stated in \cite[Proposition 6.18]{Fre2} (by replacing the $Y_{i,aq^n}$ by $Y_{i,aq^{-n}}^{-1}$, we get the $q$-character of a fundamental representation).
\qed

Now let us prove Proposition \ref{spesmalllem} :

\demo Let us denote $(1_R)$ (resp. $(2_R)$, $(3_R)$) the condition that the statement (1) (resp. (2), (3)) of the Proposition is satisfied for any $R'\leq R$. We prove by induction on $R$ simultaneously that $(1_R)$, $(2_R)$ and $(3_R)$ are satisfied. For $R = 0$ this is clear.

\noindent Now we prove the following for $R\geq 1$ :
\begin{itemize}
\item ($(1_{R-1})$ and $(2_{R-1})$ and $(3_{R-1})$) implies $(1_R)$,

\item ($(1_R)$ and $(2_{R-1})$ and $(3_{R-1})$) implies $(2_R)$,

\item ($(1_R)$ and $(2_R)$ and $(3_{R-1})$) implies $(3_R)$.
\end{itemize}
Let us start with : ($(1_R)$ and $(2_{R-1})$ and $(3_{R-1})$) implies $(2_R)$. 

\noindent By Lemma \ref{produit} 
$$\mathcal{M}(L(m))\subset (mY_{i_R,aq^{l_R}}^{-1}\mathcal{M}(V_{i_R}(aq^{l_R})))\cup (\mathcal{M}(L(mY_{i_R,aq^{l_R}}^{-1}))Y_{i_R,aq^{l_R}}).$$
As all monomials of $mY_{i_R,aq^{l_R}}^{-1}(\chi_q(V_{i_R}(aq^{l_R}))-Y_{i_R,aq^{l_R}})$ are lower than $mA_{i_R,aq^{l_R+1}}^{-1}$ (Theorem \ref{formerkr}) which is right-negative, they are not dominant. Consider $m'\in(\mathcal{M}(L(mY_{i_R,aq^{l_R}}^{-1}))Y_{i_R,aq^{l_R}}-\{m\})$. If $v_{i_R,aq^{l_R+1}}(m'm^{-1})\geq 1$ or $v_{i_R,aq^{l_R-1}}(m'm^{-1})\geq 1$, it follows from the $(1_R)$ that $m'$ is lower than $mA_{i_R,aq^{l_R+1}}^{-1}$ which is right-negative, so $m'$ is not dominant. We suppose that 
$$v_{i_R,aq^{l_R-1}}(m'm^{-1}) = v_{i_R,aq^{l_R+1}}(m'm^{-1})=0.$$ 
So we have $u_{i_R,aq^{l_R}}(m'Y_{i_R,aq^{l_R}}^{-1})\geq 0$. By $(2_{R-1})$, the monomial $m'Y_{i_R,aq^{l_R}}^{-1}\in\mathcal{M}(L(mY_{i_R,aq^{l_R}}^{-1}))$ is not dominant. So there is $i\in I$, $b\in\CC^*$, such that $(i,b)\neq (i_R,aq^{l_R})$ and $u_{i,b}(m'Y_{i_R,aq^{l_R}}^{-1}) < 0$. As $u_{i,b}(m'Y_{i_R,aq^{l_R}}^{-1})=u_{i,b}(m')$, $m'$ is not dominant. So $(2_R)$ is satisfied.

Now let us prove : ($(1_R)$ and $(2_R)$ and $(3_{R-1})$) implies $(3_R)$. 

\noindent From property $(2_R)$ and Proposition \ref{proofthin}, it suffices to prove that all monomials of $\mathcal{M}(L(m))$ are thin. Suppose that there is a monomial in $\mathcal{M}(L(m))$ which is not thin. From Lemma \ref{thinmon}, we can suppose that there is $m'\in\mathcal{M}(L(m))$ such that there are $i\in I, a\in\CC^*$ satisfying $u_{i,a}(m')=2$ and such that all $m''$ satisfying $v(m''m^{-1}) < v(m'm^{-1})$ is thin.  In particular from Proposition \ref{jdecomp} : 

$m'$ is $(\{1,\cdots,i-2\}\cup \{i\}\cup \{i+2,\dots,n\})$-dominant,

$(u_{i-1,b}(m')<0\Rightarrow b=aq)$,

$(u_{i+1,b}(m')<0\Rightarrow b=aq)$. 

\noindent (Otherwise we could construct $m''\in\mathcal{M}(L(m))$ not thin  such that $v(m''m^{-1}) < v(m'm^{-1})$). We can apply Lemma \ref{remont} for $\Glie_{\{1,\dots,i-1\}}$ of type $A_{i-1}$ and then for $\Glie_{\{i+1,\dots,n\}}$ of type $A_{n-i}$. We get a monomial $M\in\mathcal{M}(L(m))$, and by construction 

$M$ is $I-\{i\}$-dominant, 

$u_{j_1,aq^{j_1-i}}(M)\geq 1$ with $j_1\leq i$,

$u_{j_2,aq^{i-j_2}}(M)\geq 1$ with $j_1 < j_2$, $i\leq j_2$. 

\noindent Moreover as $u_{i,a}(m') = 2$, by construction $M$ is dominant. From property $(2_R)$ we have $m=M$. In particular there are $r < r'$ such that $(i_r,l_r)=(j_1,j_1-i)$ and $(i_{r'},l_{r'})=(j_2,i-j_2)$. We have 
$$l_{r'}-l_r = 2i-j_2-j_1 = i_r+i_{r'}+2(i-j_2)-2j_1 < i_r+i_{r'}.$$ 
But we have 
\begin{equation*}
\begin{split}
l_{r'}-l_r &= (l_{r'}-l_{r'-1})+\cdots +(l_{r+1}-l_r)
           \\&\geq i_{r'}+2(i_{r'-1}+\cdots + i_{r+1})+i_r\geq i_{r'}+i_r,
\end{split}
\end{equation*}
contradiction. So $(3_R)$ is satisfied.

Finally we prove : ($(1_{R-1})$ and $(2_{R-1})$ and $(3_{R-1})$) implies $(1_R)$. 

\noindent We prove $(1_R)$ by induction on $v(m'm^{-1})\geq 0$. For $v(m'm^{-1})=0$ we have $m'=m$ and the result is clear. In general consider a monomial $m' < m$ such that for $m''$ satisfying $v(m''m^{-1}) < v(m'm^{-1})$, the property $(1_R)$ is satisfied. We suppose that moreover the the property is not satisfied for $m'$, that is to say that $v_{i_R,aq^{l_R-1}}(m'm^{-1})\geq 1$ and $v_{i_R,aq^{l_R+1}}(m'm^{-1}) = 0$. It follows from Proposition \ref{jdecomp} and the induction hypothesis on $v$ that $m'$ is $(I-\{i_R\})$-dominant (otherwise we could construct $m''$ such that $v(m''m^{-1}) < v(m'm^{-1})$ and the property is not satisfied for $m''$). 

\noindent If $m'$ is not dominant, $m'$ is not $i_R$-dominant and so it follows from Proposition \ref{jdecomp} that there is $m''\in\mathcal{M}(L(m))$ $i_R$-dominant such that $m''>m'$ and $m'$ is a monomial of $L_{i_r}(m'')$. Moreover there is $b\in\CC^*$ such that $m'\leq m' A_{i_R,b} \leq m''$, and $m'A_{i_R,b}$ is a monomial of $L_{i_R}(m'')$ and so in $\mathcal{M}(L(m))$. By the induction hypothesis on $v$, $m'A_{i_R,b}$ satisfied the property $(1_R)$, and so we have $b = aq^{l_R-1}$. So $v_{i_R,aq^{l_R-1}}(m''m^{-1}) = v_{i_R,aq^{l_R+1}}(m''m^{-1}) = 0$. In particular $u_{i_R,aq^{l_R}}(m'')\geq u_{i_R,aq^{l_R}}(m)\geq 1$. By Lemma \ref{produit}, we have $m'\in \mathcal{M}(L(Y_{i_R,aq^{l_R}}))\mathcal{M}(L(m Y_{i_R,aq^{l_R}}^{-1}))$. But by Theorem \ref{formerkr} the monomials of $\mathcal{M}(L(Y_{i_R,aq^{l_R}}))$ not equal to $Y_{i_R,aq^{l_R}}$ are lower than $Y_{i_R,aq^{l_R}}A_{i_R,aq^{L_R+1}}^{-1}$. So we have $m'\in Y_{i_R,aq^{l_R}}\mathcal{M}(L(m Y_{i_R,aq^{l_R}}^{-1}))$. By the properties $(2_{R-1})$ and $(3_{R-1})$, $L(m Y_{i_R,aq^{l_R}}^{-1})$ is special and thin. In particular $u_{i_R,aq^{l_R-2}}(m'')\leq 1$, and so by Proposition \ref{aidesldeux}, $m'$ is not a monomial of $L_{i_R}(m'')$, contradiction. So $m'$ is dominant. 

\noindent As $L(mY_{i_R,aq^{l_R}}^{-1})$ is special, the monomial $m' Y_{i_R,aq^{l_R}}^{-1}$ is not dominant. So 

$u_{i_R,aq^{l_R}}(m' Y_{i_R,aq^{l_R}}^{-1}) = -1$,

$u_{j,b}(m'Y_{i_R,aq^{l_R}}^{-1}) < 0 \Rightarrow (j = i_R\text{ and }b = a q^{l_R})$. 

\noindent So we can use Lemma \ref{stara} for the thin module $L(m' Y_{i_R,aq^{l_R}}^{-1})$. Let $\alpha$, $\beta$ as in the proof of Lemma \ref{stara}. Let $j = i_R + \beta - \alpha$ and $b = aq^{l_R - \alpha - \beta -2}$. By construction of $m$ from $m'$ in the proof of Lemma \ref{stara}, we have $u_{j,b}(mY_{i_R,aq^{l_R}}^{-1})\geq 1$ and $m'\in m Y_{j,b}^{-1}\mathcal{M}(V_j(b))$. Moreover there is $R' < R$ such that $j = i_{R'}$ and $l_R-\alpha-\beta - 2 = l_{R'}$. We have 
\begin{equation*}
\begin{split}
\alpha + \beta + 2 = l_R - l_{R'} &\geq i_R+2i_{R-1}+\cdots + 2i_{R'+1}+ i_{R'}
\\&\geq 2(i_{R}+\cdots +i_{R'+1})+\beta-\alpha.
\end{split}
\end{equation*}
So $i_R+\cdots +i_{R'+1} \leq \alpha +1$ and $(i_R-\alpha)+i_{R-1}+\cdots +i_{R'+1} \leq 1$. As $i_R - \alpha \geq 1$, we have $i_{R-1} + \cdots + i_{R'+1} = 0$, $R'=R-1$ and $i_R-\alpha = 1$. 

\noindent By construction, we have $m'm^{-1}\in\ZZ[A_{i,aq^r}^{-1}]_{i\leq i_R+\beta,r\in\ZZ}$. So from Lemma \ref{depart} we can suppose that $i_R + \beta = n$. 

\noindent We have $i_R=n+1-{i_{R-1}}$. As 
\begin{equation*}
\begin{split}
\omega(m(m')^{-1})&=(\alpha_1+\cdots +\alpha_n)+(\alpha_2+\cdots +\alpha_{n-1})
\\&+\cdots +(\alpha_{i_{R-1}}+\cdots +\alpha_{n+1-{i_{R-1}}})
\\&=\alpha_1+\alpha_n+2(\alpha_2+\alpha_{n-2})+\cdots +{i_{R-1}}(\alpha_{i_{R-1}}+\alpha_{n+1-{i_{R-1}}})
\\&+{i_{R-1}}(\alpha_{{i_{R-1}}+1}+\cdots +\alpha_{n-{i_{R-1}}}),
\end{split}
\end{equation*}
the monomial $m'(mY_{{i_{R-1}},aq^{l_{R-1}}}^{-1})^{-1}$ is the lowest monomial of $\mathcal{M}(V_{i_{R-1}}(aq^{l_{R-1}}))$ (the weight of the lowest weight of fundamental representations has been computed in \cite[Lemma 6.8]{Fre2}).

\noindent Let us prove that 
\begin{equation}\label{staraidea}\mathcal{M}(L(m))\cap m'\ZZ[A_{i_R,d}]_{d\in\CC^*}\subset \{m',m'A_{i_R,aq^{l_R-1}}\}.\end{equation}
Let $m''\in(\mathcal{M}(L(m))\cap m'\ZZ[A_{i_R,d}]_{d\in\CC^*})$ different from $m'$. In particular $m\geq m'' > m'$. By construction of $m'$ from $m$, as $R'=R-1$, we have for $k \neq i_{R-1}$, $v_{k,aq^{l_k+1}}(m'm^{-1}) = 0$. So by Theorem \ref{formerkr} (for fundamental representations, that is to say the particular case proved in \cite{Fre2}), $m',m''\in m Y_{{i_{R-1}},aq^{l_{R-1}}}^{-1}\mathcal{M}(V_{i_{R-1}}(aq^{l_{R-1}}))$. As $m'(mY_{{i_{R-1}},aq^{l_{R-1}}}^{-1})^{-1}$ is the lowest monomial of $\mathcal{M}(V_{i_{R-1}}(aq^{l_{R-1}}))$, Lemma \ref{fundalgo} gives :
\begin{equation*}
\begin{split}
(\chi_q(V_{i_{R-1}}(aq^{l_{R-1}}))Y_{i_R,aq^{L_R}} -1)A_{i_R,aq^{l_R-1}}^{-1} \in& 1+A_{i_R+1,aq^{l_R-2}}\ZZ[A_{k,d}]_{k\in I,d\in\CC^*}
\\&+A_{i_R-1,aq^{l_R-2}}\ZZ[A_{k,d}]_{k\in I,d\in\CC^*}.
\end{split}
\end{equation*}
As by hypothesis $v_{i_R - 1}(m'(m'')^{-1}) = v_{i_R + 1}(m'(m'')^{-1}) = 0$, we get :
$$m''(mY_{{i_{R-1}},aq^{l_{R-1}}}^{-1})^{-1} = m'A_{i_R,aq^{l_R-1}}(mY_{{i_{R-1}},aq^{l_{R-1}}}^{-1})^{-1}.$$
Let us prove that 
\begin{equation}\label{staraideb} m''\in (\mathcal{M}(L(m))\cap m'\ZZ[A_{i_R,d}^{\pm 1}]_{d\in\CC^*})-\{m'A_{i_R,b}\}\Rightarrow v_{i_R}(m''(m')^{-1})\geq 0
.\end{equation}
Consider a monomial $m''$ satisfying the left property of (\ref{staraideb}). By Lemma \ref{plusgrand}, for $k\neq R-1$ we have $v_{i_k,aq^{l_k+1}}(m'm^{-1}) = 0$. So 
$$m''\in m (Y_{i_{R-1},aq^{l_{R-1}}}^{-1}\mathcal{M}(V_{i_{R-1}}(aq^{l_{R-1}})))(\prod_{\{k|i_k=i_R\}}Y_{i_R,aq^{l_k}})^{-1}\prod_{\{k|i_k=i_R\}}\mathcal{M}(V_{i_R}(aq^{l_k})).$$ 
Let us write this decomposition $$m''=mY_{i_{R-1},aq^{l_{R-1}}}^{-1}(m'')_{R-1}(\prod_{\{k|i_k=i_R\}}Y_{i_R,aq^{l_k}})^{-1}\prod_{\{k|i_k=i_R\}}(m'')_k.$$ 
(If $i_{R-1} = i_R$ we put $(m'')_{R-1}$ only one time). Let $k\neq R-1$ satisfying
$i_k=i_R$. Observe that for $R_1<R_2$, we have $l_{R_2}-l_{R_1}\geq i_{R_1}+i_{R_2}\geq 2$. So by Lemma \ref{plusgrand}, $v_{i_R+1,aq^{l_k+2}}(m'm^{-1})=v_{i_R - 1,aq^{l_k+2}}(m'm^{-1})=0$. So $(m'')_k = Y_{i_R,aq^{l_k}}$ or $(m'')_k=Y_{i_R,aq^{l_k}}A_{i_R,aq^{l_k+1}}^{-1}$. As a consequence, $(m'')_{R-1} = Y_{i_R,aq^{l_R}}^{-1}$ or $(m'')_{R-1}=Y_{i_R,aq^{l_R}}^{-1}A_{i_R,aq^{l_R-1}}$ (Lemma \ref{plusgrand}). So
\begin{equation*}
\begin{split}
v_{i_R}(m''(m')^{-1}) &= v_{i_R}((m'')_{R-1}Y_{i_R,aq^{l_R}}) + \sum_{\{k\neq R-1|i_k=i_R\}}v_{i_R}((m'')_k Y_{i_R,aq^{l_k}}^{-1})
\\&\geq v_{i_R}((m'')_{R-1}Y_{i_R,aq^{l_R}}) \geq -1.
\end{split}
\end{equation*}
If $v_{i_R}(m''(m')^{-1})=-1$, then for all $k$ satisfying $i_k=i_R$ we have $(m'')_k=Y_{i_R,aq^{l_k}}$ and $(m'')_{R-1}=Y_{i_R,aq^{l_R}}^{-1}A_{i_R,aq^{l_R-1}}$. So $m''=m'A_{i_R,aq^{l_{R-1}}}$ and we can conclude the proof of (\ref{staraideb}).

\noindent Now to prove it suffices to prove that the conditions of Theorem \ref{racourc} with $i = i_R$ are satisfied for $m'$. 

Condition (i) of Theorem \ref{racourc} : 
\\the unicity follows from the statement (\ref{staraidea}) above. For the existence, it suffices to prove that $M = m'A_{i_R,aq^{l_R-1}}$ is in $\mathcal{M}(L(m))$. By Lemma \ref{oudeux}, there is $j\in I$, $M'\in\mathcal{M}(L(m))$ $j$-dominant such that $M'>m'$ and $M'\in m'\ZZ[A_{j,a}]_{a\in\CC^*}$. By the induction hypothesis on $v$ we have $j = i_R$, and so by unicity $M'=M$. 

Condition (ii) of Theorem \ref{racourc} :
\\we have by Lemma \ref{ouun}
$$\sum_{r\in\ZZ} x_{i,r}^+(V_M) \subset \sum_{m'\in m\ZZ[A_{i,d}^{\pm}]_{d\in\CC^*}}(L(m))_{m'},$$ 
and so the result follows from the statement (\ref{staraideb}) above.

Condition (iii) of Theorem \ref{racourc} : 
\\by Lemma \ref{produit}, we have $M\in \mathcal{M}(L(Y_{i_R,aq^{l_R}}))\mathcal{M}(L(m Y_{i_R,aq^{l_R}}^{-1}))$. But by Theorem \ref{formerkr} the monomials of $\mathcal{M}(L(Y_{i_R,aq^{l_R}}))$ not equal to $Y_{i_R,aq^{l_R}}$ are lower than $Y_{i_R,aq^{l_R}}A_{i_R,aq^{L_R+1}}^{-1}$. So we have $M\in Y_{i_R,aq^{l_R}}\mathcal{M}(L(m Y_{i_R,aq^{l_R}}^{-1}))$. By $(3_{R-1})$, the module $L(mY_{i_R,aq^{l_R}}^{-1})$ is thin and so $u_{i_R,aq^{l_R-2}}(M)\leq 1$. Moreover by the induction hypothesis on $v$, $v_{i_R,aq^{l_R-1}}(Mm^{-1})=v_{i_R,aq^{l_R+1}}(Mm^{-1})=0$. So $u_{i_R,aq^{l_R}}(M)\geq 1$. So by Proposition \ref{aidesldeux}, $m'$ is not a monomial of $\mathcal{M}(L_{i_R}(M))$. 

Condition (iv) of Theorem \ref{racourc} : 
\\the result follows from the statement (\ref{staraideb}) above.

Condition (v) of Theorem \ref{racourc} : 
\\clear by the induction hypothesis on $v$.

\qed

The case of standard modules $M(X_{k,a}^{(n)})$ can be studied in the same way by replacing $i$ by $\overline{i}=n-i+1$.

We can conclude the proof of Proposition \ref{equismalla} with Proposition \ref{fundsmall}, Proposition \ref{kegaldeux} and the following counter examples :

\begin{lem}\label{reste} We suppose that $n\geq 3$. Let $k\geq 3, a\in\CC^*$ and $1<i<n$. Then $M(X_{k,a}^{(i)})$ is not small.\end{lem}

\demo Consider $m'=X_{k,a}^{(i)}A_{i,aq^{-2k+2}}^{-1}\leq X_{k,a}^{(i)}$. Then $m'$ is dominant. As $\Glie_{\{i-1 ,i ,i+1\}}$ is of type $sl_4$, by using Lemma \ref{depart}, we can check as in remark \ref{counter} that $L(m')$ is not special, and so $M(X_{k,a}^{(i)})$ is not small.\qed

\subsection{End of the proof of Theorem \ref{mainres}}\label{odiss} In general for $\Glie$ not of type $A$, $i$ extremal does not imply that $M(X_{k,a}^{(i)})$ is small. For example :

\begin{rem}\label{counterd} Let $\Glie$ be of type $D_4$ and $m=Y_{1,q^3}Y_{1,q^5}Y_{2,1}$. By using the process described in remark \ref{process}, the following monomials occur in $\chi_q(L(m))$ : $1_3 1_5 2_1$, $1_1 1_3 1_5 2_2^{-1} 3_1 4_1$, $1_1 1_3 1_5 2_2 3_3^{-1} 4_3^{-1}$, $1_1 1_3^2 1_5 2_4^{-1}$, $1_1 1_3$. So $L(m)$ is not special. As $Y_{1,q^3}Y_{1,q^5}Y_{2,1}=X_{4,q^2}^{(1)}A_{1,1}^{-1}\in \mathcal{M}(M(X_{4,q^2}^{(1)}))$, $M(X_{4,q^2}^{(1)})$ is not small. \end{rem}

Let us end the proof of Theorem \ref{mainres} :

The case $k = 1$ follows from Lemma \ref{fundsmall}. The case $k = 2$ follows from Lemma \ref{kegaldeux}. In the rest of the proof we suppose that $k\geq 3$. 

\noindent Suppose that $i$ is not extremal. There are $j\neq j'$ such that $C_{i,j}=C_{i,j'}=-1$. Consider $m' = X_{k,a}^{(i)}A_{i,aq^{-2k+2}}^{-1}\leq X_{k,a}^{(i)}$. Then $m'$ is dominant. Let $J=\{i,j,j'\}$. $\Glie_J$ is of type $A_3$ and so by using Lemma \ref{depart}, we can check as in remark \ref{counter} that $L(m')$ is not special.

\noindent Suppose that $i$ is extremal. Let $i_2$ be the unique element of $I$ satisfying $C_{i,i_2} = -1$. Let $i_3,\cdots, i_{d_i}$ such that for $2\leq r\leq d_i-1$, $C_{i_r,i_{r+1}}=-1$ and $i_{d_i}$ is special. Let $i_{d_i+1}\neq i_{d_i+2}$ such that $C_{i_{d_i},i_{d_i+1}}=C_{i_{d_i},i_{d_i+2}}=-1$ and $i_{d_i-1}, i_{d_i+1}, i_{d_i+2}$ are distinct.

\noindent For illustration an example is given on the following picture :

{\large 
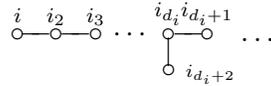
\vspace{-.15cm} 
$$\stackrel{i}{\circ}\hspace{-.18cm}\sn\hspace{-.18cm}\stackrel{i_2}{\circ} \hspace{-.18cm}\sn\hspace{-.18cm}\stackrel{i_3}{\circ} 
\dots\stackrel{i_{d_i}}{\circ}\hspace{-.24cm}\sn\hspace{-.42cm}\stackrel{i_{d_i + 1}}{\circ}
\dots$$
\vspace{-.85cm}$$\hspace{.62cm}|$$
\vspace{-.85cm}$$\hspace{1.42cm}\circ\text{ \scriptsize $i_{d_i + 2}$}$$}

\noindent Suppose that $k\geq d_i + 2$. Let $m' = X_{k,a}^{(i)}A_{i,aq^{2-k}}^{-1} = Y_{i_2,aq^{2-k}}X_{k-1,aq}^{(i)}$. By remark \ref{process}, 
\begin{equation*}
\begin{split}
m'' = &m'(A_{i_2,aq^{3-k}}^{-1}A_{i_3,aq^{4-k}}^{-1} \cdots A_{i_{d_i+1},aq^{d_i+2-k}}^{-1})
\\&\times (A_{i_{d_i+2},aq^{d_i+2-k}}^{-1}A_{i_{d_i},aq^{d_i+3-k}}^{-1}A_{i_{d_i-1},aq^{d_i+4-k}}^{-1}\cdots A_{i_2,aq^{2d_i+1-k}}^{-1})\in\mathcal{M}(L(m')).
\end{split}
\end{equation*}
But 
\begin{equation*}
\begin{split}
(m'')^{\rightarrow(i)} &= X_{k-1,aq}^{(i)}(A_{i_2,aq^{3-k}}^{-1}A_{i_2,aq^{2d_i+1-k}}^{-1})^{\rightarrow(i)} 
\\&= Y_{i,aq^{3-k}}Y_{i,aq^{5-k}}\cdots Y_{i,aq^{k - 1}}Y_{i,aq^{2d_i+1-k}},
\end{split}
\end{equation*}
and $3-k\leq 2d_i + 1 - k\leq k - 3$. So $m''A_{i,aq^{2d_i+2-k}}^{-1}$ is dominant and occurs in $\chi_q(L(m''))$. So $L(m'')$ is not special and $M(X_{k,a}^{(i)})$ is not small.

\noindent Suppose that $k \leq d_i + 1$ and there there is a dominant monomial $m' < X_{k,a}^{(i)}$. By Lemma \ref{lequeldom}, $(v_j(m'm^{-1}) \neq 0 \Rightarrow j\in\{i_1,\cdots,i_{d_i}\})$. So from Lemma \ref{depart}, we can work with $\Glie_{\{i_1,\cdots, i_{d_i}\}}$ of type $A_{d_i}$. So it follows from Proposition \ref{equismalla} that $M(X_{k,a}^{(i)})$ is small.\qed

\subsection{General simply laced quantum affinizations}

The notion of quantum affinization can be extended beyond quantum affine algebras : the quantum affinization $\U_q(\hat{\Glie})$ of a quantum Kac-Moody algebra $\U_q(\Glie)$ is defined with the same generators and relations as the Drinfeld realization of quantum affine algebras, but by using the generalized symmetrizable Cartan matrix of $\Glie$ instead of a Cartan matrix of finite type. The quantum affine algebra, quantum affinizations of usual quantum groups, are the simplest examples have the particular property of being also quantum Kac-Moody algebras. The quantum affinization of a quantum affine algebra is called a quantum toroidal algebra (or double affine quantum algebra). It is not a quantum Kac-Moody algebra, but is also of particular interest, in particular in relation to double affine Hecke algebras (Cherednik algebras).

In \cite{mi, Naams, her04}, the category $\mathcal{O}$ of integrable representations is studied. One can define for general quantum affinizations analogs of Kirillov-Reshetikhin modules (these representations are not finite dimensional in general). We can also define the notion of small modules by using the characterization in Theorem \ref{repthchar}.

The statement of Theorem \ref{mainres} is satisfied for all simply-laced quantum affinizations, by using exactly the same proof, except that in the end of the proof of Theorem \ref{mainres} (subsection \ref{odiss}), for $J=\{i,j,j'\}$, $\Glie_J$ may be of type $A_3$ or of type $A_2^{(1)}$ (in the second case we have $C_{i,j} = C_{i,j'} = C_{j,j'} = -1$). In this case and we can check as in the following remark that for $m'$ as in subsection \ref{odiss}, $L(m')$ is not small. 

\begin{rem}\label{counterloop} Let $\Glie$ be of $A_2^{(1)}$, consider $m=Y_{2,1}Y_{2,q^2}Y_{2,q^4}$, $m' = mA_{2,q}^{-1} = Y_{1,q}Y_{0,q}Y_{2,q^4}$. Then by using the process described in remark \ref{process}, the monomials $$Y_{1,q^3}^{-1}Y_{0,q^3}^{-1}Y_{1,q^2}Y_{0,q^2}Y_{2,q^2}^2Y_{2,q^4}=m'A_{1,q^2}^{-1}A_{0,q^2}^{-1},$$ $$Y_{2,q^2}Y_{1,q^2}Y_{0,q^2}=m'A_{1,q^2}^{-1}A_{0,q^2}^{-1}A_{2,q^3}^{-1},$$ 
occur in $\chi_q(L(m'))$ and $L(m')$ is not special. So $M(m)$ is not small.\end{rem}


\begin{thebibliography}{99}

\bibitem[B]{bou} {\bf N. Bourbaki}, {\it Groupes et alg\`ebres de Lie}, {Chapitres IV-VI, Hermann (1968)}

\mk

\bibitem[BBD]{bbdg} {\bf A. Beilinson, J. Bernstein and P. Deligne}, {\it Faisceaux pervers}, {Analysis and topology on singular spaces, I (Luminy, 1981),  5--171, Astérisque, {\bf 100}, Soc. Math. France, Paris (1982)}

\mk

\bibitem[BM]{bm} {\bf W. Borho and R. MacPherson}, {\it Partial resolutions of nilpotent varieties}, {Analysis and topology on singular spaces, II, III (Luminy, 1981),  23--74, Astérisque, {\bf 101-102}, Soc. Math. France, Paris (1983)}

\mk

\bibitem[CM1]{cm} {\bf V. Chari and A. Moura}, {\it Characters and blocks for finite-dimensional representations of quantum affine algebras}, {Int. Math. Res. Not. {\bf 2005}, no. 5, 257--298 (2005)}

\mk

\bibitem[CM2]{cm2} {\bf V. Chari and A. Moura}, {\it Characters of fundamental representations of quantum affine algebras}, {Acta Appl. Math. {\bf 90},  no. 1-2, 43--63 (2006)}

\mk

\bibitem[CP1]{Cha0} {\bf V. Chari and A. Pressley}, {\it Quantum Affine Algebras}, {Comm. Math. Phys. {\bf 142}, 261-283 (1991)}

\mk

\bibitem[CP2]{Cha}{\bf V. Chari and A. Pressley}, {\it Quantum affine algebras and their representations}, {in Representations of groups (Banff, AB, 1994),59-78, CMS Conf. Proc, {\bf 16}, Amer. Math. Soc., Providence, RI (1995)}

\mk

\bibitem[CP6]{Cha2}{\bf V. Chari and A. Pressley}, {\it A Guide to Quantum Groups}, {Cambridge University Press, Cambridge (1994)} 

\mk

\bibitem[CP7]{Cha4}{\bf V. Chari and A. Pressley}, {\it Integrable and Weyl modules for quantum affine ${\rm sl}\sb 2$}, {Quantum groups and Lie theory (Durham, 1999), 48--62, London Math. Soc. Lecture Note Ser., {\bf 290}, Cambridge Univ. Press, Cambridge, (2001)}

\mk

\bibitem[Dr1]{Dri1}{\bf V. G. Drinfeld}, {\it Quantum groups}, {Proceedings of the International Congress of Mathematicians, Vol. 1, 2 (Berkeley, Calif., 1986), 798--820, Amer. Math. Soc., Providence, RI, (1987)}

\mk

\bibitem[Dr2]{Dri2}{\bf V. G. Drinfeld}, {\it A new realization of Yangians and of quantum affine algebras}, {Soviet Math. Dokl. {\bf 36}, no. 2, 212--216 (1988)}

\mk

\bibitem[DM]{dm}{\bf G. W. Delius and N. J. MacKay}, {\it Affine quantum groups}, {Encyclopedia of Mathematical Physics, Elsevier (2006)}

\mk

\bibitem[FM]{Fre2} {\bf E. Frenkel and E. Mukhin}, 
{\it Combinatorics of $q$-Characters of Finite-Dimensional Representations of Quantum Affine Algebras}, {Comm. Math. Phy., vol {\bf 216}, no. 1, pp 23-57 (2001)}

\mk

\bibitem[FR]{Fre} {\bf E. Frenkel and N. Reshetikhin}, {\it The $q$-Characters of Representations of Quantum Affine Algebras and Deformations of $W$-Algebras}, {Recent Developments in Quantum Affine Algebras and related topics, Cont. Math., vol. {\bf 248}, 163--205 (1999)}

\mk

\bibitem[GM1]{gm1} {\bf M. Goresky and R. MacPherson}, {\it Intersection homology theory}, {Topology {\bf 19}, no. 2, 135--162 (1980)}

\mk

\bibitem[GM2]{gm2} {\bf M. Goresky and R. MacPherson}, {\it Intersection homology. II}, {Invent. Math. {\bf 72}, no. 1, 77--129 (1983)} 

\mk

\bibitem[Ha]{ha} {\bf R. Hardt}, {\it Semi-algebraic local-triviality in semi-algebraic mappings}, {Amer. J. Math. {\bf 102}, no. 2, 291--302 (1980)}

\mk

\bibitem[He1]{her02} {\bf D. Hernandez}, {\it Algebraic approach to $q,t$-characters}, {Adv. Math.  {\bf 187},  no. 1, 1--52 (2004)}

\mk

\bibitem[He2]{her04} {\bf D. Hernandez}, {\it Representations of quantum affinizations and fusion product}, {Transform. Groups {\bf 10}, no. 2, 163--200 (2005)}

\mk

\bibitem[He3]{her05} {\bf D. Hernandez}, {\it Monomials of q and q,t-chraracters for non simply-laced quantum 
affinizations}, {Math. Z. {\bf 250}, no. 2, 443--473 (2005)}

\mk

\bibitem[He4]{her06} {\bf D. Hernandez}, {\it The Kirillov-Reshetikhin conjecture and solutions of T-systems}, {J. Reine Angew. Math. {\bf 596}, 63--87 (2006)}

\mk

\bibitem[He5]{her07} {\bf D. Hernandez}, {\it Drinfeld coproduct, quantum fusion tensor category and applications}, {Proc. London Math. Soc. (3) {\bf 95}, no. 3, 567--608 (2007)}

\mk

\bibitem[He6]{miniaff} {\bf D. Hernandez}, {\it On minimimal affinizations of representations of quantum groups}, {Comm. Math. Phys. {\bf 277}, no. 1, 221--259 (2007)}

\mk

\bibitem[J]{jim} {\bf M. Jimbo}, {\it A $q$-difference analogue of $\U({\Glie})$ and the Yang-Baxter equation}, {Lett. Math. Phys. {\bf 10}, no. 1, 63--69 (1985)}

\mk

\bibitem[Ka]{kac} {\bf V. Kac}, {\it Infinite dimensional Lie algebras}, {3rd Edition, Cambridge University Press (1990)}

\mk

\bibitem[Kn]{kn} {\bf H. Knight}, {\it Spectra of tensor products of finite-dimensional representations of Yangians}, {J. Algebra {\bf 174}, no. 1, 187--196 (1995)}

\mk

\bibitem[KR]{kr} {\bf A.N. Kirillov and N. Reshetikhin}, {\it Representations of Yangians and multiplicities of the inclusion of the irreducible components of the tensor product of representations of simple Lie algebras}, {J. Soviet Math. {\bf 52}, no. 3, 3156--3164 (1990); translated from Zap. Nauchn. Sem. Leningrad. Otdel. Mat. Inst. Steklov. (LOMI) 160, Anal. Teor. Chisel i Teor. Funktsii. 8, 211--221, 301 (1987)}

\mk

\bibitem[KS]{ks} {\bf A. Kuniba and S. Suzuki}, {\it Analytic Bethe Ansatz for fundamental representations and yangians}, {Commun. Math. Phys. {\bf 173}, 225 - 264 (1995)}

\mk

\bibitem[M]{mi} {\bf K. Miki}, {\it Representations of quantum toroidal algebra $U\sb q({\rm sl}\sb {n+1,{\rm tor}}) (n\geq 2)$}, {J. Math. Phys. {\bf 41}, no. 10, 7079--7098 (2000)}

\mk

\bibitem[N1]{nun} {\bf H. Nakajima}, {\it Instantons on ALE spaces, quiver varieties, and Kac-Moody algebras}, {Duke Math. J. {\bf 76},  no. 2, 365--416 (1994)}

\mk

\bibitem[N2]{NaDuke} {\bf H. Nakajima}, {\it Quiver varieties and Kac-Moody algebras}, {Duke Math. J. {\bf 91}, no.~3, 515--560 (1998)}.

\mk

\bibitem[N3]{Naams} {\bf H. Nakajima}, {\it Quiver varieties and finite-dimensional representations of quantum affine algebras}, {J. Amer. Math. Soc. {\bf 14},  no. 1 (2001)}

\mk

\bibitem[N4]{Naex} {\bf H. Nakajima}, {\it $T$-analogue of the $q$-characters of finite dimensional representations of quantum affine algebras}, {Physics and combinatorics, 2000 (Nagoya),  196--219, World Sci. Publishing, River Edge, NJ (2001)}

\mk

\bibitem[N5]{NaCong} {\bf H. Nakajima}, {\it Geometric construction of representations of affine algebras}, {Proceedings of the International Congress of Mathematicians, Vol. I (Beijing, 2002),  423--438, Higher Ed. Press, Beijing (2002)}

\mk

\bibitem[N6]{Nae}{\bf H. Nakajima}, {\it Problems on quiver varieties}, {The 50th Geometry Symposium, Hokkaido Univ. 2003 Aug., online version http://www.math.kyoto-u.ac.jp/\textasciitilde nakajima/TeX/kika03.pdf}

\mk

\bibitem[N7]{Nad}{\bf H. Nakajima}, {\it $t$-analogs of $q$-characters of Kirillov-Reshetikhin modules of quantum 
affine algebras}, {Represent. Theory {\bf 7}, 259--274 (electronic) (2003)}

\mk

\bibitem[N8]{Nab} {\bf H. Nakajima}, {\it Quiver Varieties and $t$-Analogs of $q$-Characters of Quantum Affine Algebras}, {Ann. of Math. {\bf 160}, 1057 - 1097 (2004)} 

\mk

\bibitem[S]{os} {\bf O. Schiffmann}, {\it Vari\'et\'es carquois de Nakajima}, {S\'eminaire Nicolas Bourbaki no. {\bf 976}, vol. 2006-2007 (2007)}

\mk

\bibitem[T]{t} {\bf R. Thom}, {\it Ensembles et morphismes stratifiés}, {Bull. Amer. Math. Soc. {\bf 75}, 240--284 (1969)}

\mk

\bibitem[VV]{vv} {\bf M. Varagnolo and E. Vasserot}, {\it Standard modules of quantum affine algebras}, {Duke Math. J. {\bf 111}, no. 3, 509--533 (2002)}

\end{thebibliography}
\end{document}